\documentclass{article}
\usepackage{setspace}
\usepackage[utf8]{inputenc}
\usepackage[english]{babel}
\usepackage{amsfonts,amsmath,amssymb,mathtools}
\usepackage[left=2cm,right=2cm,top=2cm,bottom=2cm,bindingoffset=0cm]{geometry}
\usepackage{graphicx}
\usepackage{authblk}
\usepackage{xcolor}
\usepackage{hyperref}
\usepackage{dsfont}
\usepackage{authblk}
\usepackage{slashed}
\usepackage{braket}
\usepackage{comment}
\usepackage{pb-diagram}
\usepackage{wasysym}
\usepackage{amsthm}
\usepackage{cancel}
\usepackage{tikz}\usetikzlibrary{tikzmark}
\usetikzlibrary{fadings} 
\usetikzlibrary{positioning}
\usetikzlibrary{backgrounds} 
\usetikzlibrary{decorations.pathreplacing}
\usetikzlibrary{arrows}
\definecolor{airforceblue}{rgb}{0.36, 0.54, 0.66}	
\definecolor{beige}{rgb}{0.96, 0.96, 0.86}
\definecolor{bittersweet}{rgb}{1.0, 0.44, 0.37}
\definecolor{melon}{rgb}{0.99, 0.74, 0.71}
\definecolor{mustard}{rgb}{1.0, 0.86, 0.35}
\definecolor{lava}{rgb}{0.81, 0.06, 0.13}
\definecolor{magnolia}{rgb}{0.97, 0.96, 1.0}
\definecolor{lavendermist}{rgb}{0.9, 0.9, 0.98}
\definecolor{lavendergray}{rgb}{0.77, 0.76, 0.82}
\definecolor{palepink}{rgb}{0.98, 0.85, 0.87}
\definecolor{palesilver}{rgb}{0.79, 0.75, 0.73}
\definecolor{cadetgrey}{rgb}{0.57, 0.64, 0.69}
\definecolor{anti-flashwhite}{rgb}{0.95, 0.95, 0.96}
\colorlet{Light0anti-flashwhite}{anti-flashwhite!70!white}
\colorlet{Lightanti-flashwhite}{anti-flashwhite!50!white}
\colorlet{Light2anti-flashwhite}{anti-flashwhite!30!white}
\definecolor{linkcolor}{rgb}{0,0,1}
\definecolor{urlcolor}{rgb}{0,0,1}

\usepackage{tikz-cd}
\usetikzlibrary{decorations.markings}
\usetikzlibrary{positioning}
\usepackage{braids}
\usepackage{array}
\usepackage{ytableau}
\usepackage{longtable}
\usepackage{empheq}
\usepackage{multicol}
\usepackage{mathrsfs}
\usepackage{diagbox}
\usepackage[vcentermath]{youngtab}
\usepackage[natbib=true, backend = bibtex, style=numeric-comp, sorting=none]{biblatex}
\hypersetup{pdfstartview=FitH,  linkcolor=linkcolor,urlcolor=urlcolor,citecolor=urlcolor, colorlinks=true}

\newcommand\bem{\begin{pmatrix}}
\newcommand\eem{\end{pmatrix}}
\newcommand\beq{\begin{equation}}
\newcommand\eeq{\end{equation}}
\newcommand\beqs{\begin{equation*}}
\newcommand\eeqs{\end{equation*}}

\newcommand{\tr}{\text{tr}}

\date{}

\def\be{\begin{eqnarray}}
\def\ee{\end{eqnarray}}

\def\tr{{\rm tr}\,}

\definecolor{red}{rgb}{1,0,0}
\definecolor{orange}{rgb}{1,0.5,0}
\definecolor{violet}{rgb}{0.7,0,1}

\newtheorem{theorem}{Theorem}[section]
\newtheorem{definition}{Definition}[section]

\newtheorem{proposition}{Proposition}[section]
\newtheorem{remark}{Remark}[section]
\newtheorem{corollary}{Corollary}[section]

\begin{document}

\title{\bf Closed 4-braids and the Jones unknot conjecture
}

\author[1]{{\bf D. Korzun}\thanks{\href{mailto:korzun.dv@phystech.edu}{korzun.dv@phystech.edu}}}
\author[1,3]{{\bf E. Lanina}\thanks{\href{mailto:lanina.en@phystech.edu}{lanina.en@phystech.edu}}}
\author[1,2,3]{{\bf A. Sleptsov}\thanks{\href{mailto:sleptsov@itep.ru}{sleptsov@itep.ru}}}

\vspace{5cm}

\affil[1]{Moscow Institute of Physics and Technology, 141700, Dolgoprudny, Russia}
\affil[2]{Institute for Information Transmission Problems, 127051, Moscow, Russia}
\affil[3]{NRC "Kurchatov Institute", 123182, Moscow, Russia\footnote{former Institute for Theoretical and Experimental Physics, 117218, Moscow, Russia}}
\renewcommand\Affilfont{\itshape\small}

\maketitle

\vspace{-7cm}

\begin{center}
	\hfill MIPT/TH-22/23\\
	\hfill ITEP/TH-29/23\\
	\hfill IITP/TH-23/23
\end{center}

\vspace{4.5cm}

\begin{abstract}

{
The Jones problem is a question whether there is a non-trivial knot with the trivial Jones polynomial in one variable $q$. The answer to this fundamental question is still unknown despite numerous attempts to explore it. In braid presentation the case of 4-strand braids is already open. S. Bigelow showed in 2000 that if the Burau representation for four-strand braids is unfaithful, then there is an infinite number of non-trivial knots with the trivial two-variable HOMFLY-PT polynomial and hence, with the trivial Jones polynomial, since it is obtained from the HOMFLY-PT polynomial by the specialisation of one of the variables $A=q^2$.

In this paper, we study four-strand braids and ask whether there are non-trivial knots with the trivial Jones polynomial but a non-trivial HOMFLY-PT polynomial. We have discovered that there is a whole 1-parameter family, parameterised by the writhe number, of 2-variable polynomials that can be HOMFLY-PT polynomials of some knots. We explore various properties of the obtained hypothetical HOMFLY-PT polynomials and suggest several checks to test these formulas. A generalisation is also proposed for the case of a large number of strands.

}
\end{abstract}

\tableofcontents

\section{Introduction}

In this paper, we apply approaches for the study of quantum knot invariants to the celebrated {\it Jones unknot conjecture}. It poses a question whether a non-trivial knot with the unit Jones polynomial exists. At first, let us discuss some motivation and current developments in this area.

A polynomial knot invariant is a function of a knot taking values in polynomials which is invariant under Reidemeister moves and ambient isotopies. One of the most famous and strongest polynomial knot invariant -- the Jones quantum polynomial $J(q)$, was introduced by Vaughan Jones in 1984~\cite{jones2005jones,jones1987hecke,jones1985polynomial}, when he studied von Neumann operators acting on a state space of some system. The Jones polynomial can be defined as a polynomial satisfying the following skein relations
\be\label{JonesSkein}
\begin{picture}(300,55)(-70,-59)

\put(0,-30){
\put(-42,10){\mbox{$q^2\cdot J\Bigg( $}}
\put(15,10){\mbox{$\Bigg) $}}
\put(-10,0){\vector(1,1){24}}
\put(14,0){\line(-1,1){10}}
\put(0,14){\vector(-1,1){10}}

\put(64,0){
\put(-38,10){\mbox{$- \ \, q^{-2}\cdot J\Bigg($}}
\put(39,10){\mbox{$\Bigg)$}}
\put(37,0){\vector(-1,1){24}}
\put(13,0){\line(1,1){10}}
\put(27,14){\vector(1,1){10}}
}

\put(152,0){
\put(-39,10){\mbox{$= \, (q-q^{-1})\cdot J\Bigg($}}
\put(53,10){\mbox{$\Bigg)\,,$}}
\put(39,0){\vector(0,1){24}}
\put(49,0){\vector(0,1){24}}
}
}

\put(3,-53){\circle{21}}
\put(-22,-56){\mbox{$J\bigg( $}}
\put(15,-56){\mbox{$\bigg) $}}
\put(27,-56){\mbox{$= \, 1\,.$}}

\end{picture}
\ee


It was a breakthrough when quantum topology and in particular the Jones polynomial turned out to be connected with various topics of theoretical and mathematical physics such as quantum field theory~\cite{Witten1988hf,Chern1974ft,Guadagnini:1989kr,GUADAGNINI1990575,Kaul1991np,RamaDevi1992np,Ramadevi1993np,Ramadevi1994zb}, quantum groups~\cite{TURAEV1992865,Reshetikhin:1990pr}, topological string theory~\cite{Ooguri1999bv, Labastida2000zp, Labastida2000yw, Labastida2001ts, Marino2001re} and conformal field theory~\cite{Kaul1991np,RamaDevi1992np,Ramadevi1993np,Zodinmawia2011oya,Zodinmawia2012kn}. That is why quantum knot invariants attract lots of attention among both physicists and mathematicians. 

\begin{figure}[h!]
		\centering	
		\includegraphics[width =0.6\linewidth]{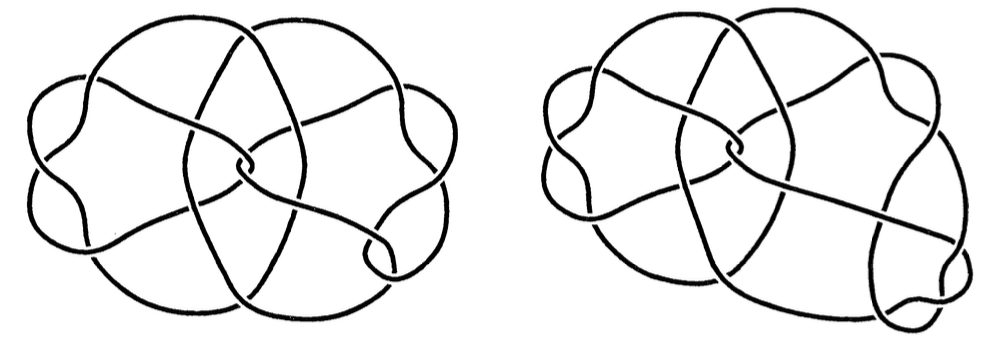}
		\caption{Two-component links with 15 and 17 crossings with the trivial Jones polynomial from~\cite{thistlethwaite2001links}.}
		\label{pic:LinksTrivJones}
	\end{figure}

The Jones polynomial does not distinguish all knots. Already in the 80s, infinite number of pairs of knots with equal Jones polynomials was revealed~\cite{birman1985jones,lickorish1987linear,kanenobu1986examples,lickorish1987polynomials}. For example, it takes the same value on knots related by mutation (the so-called mutant knots)~\cite{conway1970enumeration,morton1996distinguishing}, such as the Conway and the Kinoshita–Terasaka knots~\cite{tillmann2000kinoshita}. But a more basic question is still open. Does the Jones polynomial $J(q)$ detect the unknot? This is the famous Jones unknot problem, that is a fundamental problem on the invariant topology. Probably, absence of an answer to the Jones problem interferes with understanding the topological meaning of the variable $q$~\cite{jones2005jones}. At the same time, it is known that non-trivial links with the Jones polynomial equal to the one of unlink\footnote{We also call such Jones polynomial as the trivial Jones polynomial.} do exist~\cite{de1992adequate,ganzell2019unoriented,eliahou2003infinite,thistlethwaite2001links} (for example, see Fig.~\ref{pic:LinksTrivJones}). It is also interesting that it was only recently proved that a particular $t$-deformation of the Jones polynomial -- the Khovanov polynomial~\cite{khovanov2000categorification,bar2002khovanov}, does detect the unknot~\cite{kronheimer2011khovanov}.     
 
 Since V. Jones stated the unknot conjecture that the Jones polynomial detects the unknot, this problem has become classical. There are lots of attempts to solve the Jones problem. V. Jones and D. Rolfsen studied this problem regarding 4-braids~\cite{jones1994theorem}. They presented a method of producing a non-trivial knot knot with the trivial Jones polynomial, but did not succeed. A bit later T. Stanford proved the following theorem~\cite{stanford1996braid,stanford1998vassiliev}. Two knots have same all Vassiliev invariants of order less than $n$ if and only if they are equivalent modulo the $n-$th group of the lower central series of some pure braid group\footnote{The $n-$th group of the lower central series of the pure braid group $P_k$ is a group ${\rm LCS}_n(P_k)=[{\rm LCS}_{n-1}(P_k),P_k]$ with ${\rm LCS}_1(P_k)=P_k$. Remind that in a pure braid, the beginning and the end of each strand are in the same position.}. This statement opens a way to find a non-trivial knot with the unit Jones polynomial from some diagram of the unknot. 
 
 Computer checks were carried out for knots with up to 24 crossings~\cite{tuzun2021verification}, and no non-trivial knots with the unit Jones polynomial was found. V. Manturov put forward two hypothesis to forbid non-trivial knots with the trivial Jones polynomial by connection of the Jones problem with virtual knots~\cite{manturov2018knot}. Similar problem concerns the HOMFLY-PT polynomial\footnote{The HOMFLY-PT polynomial $H(q,A)$ can be defined with the use of the following skein relations:
\be\label{Hrelation}
\begin{picture}(300,25)(-70,-30)

\put(0,-30){
\put(-40,10){\mbox{$A\cdot H\Bigg( $}}
\put(15,10){\mbox{$\Bigg) $}}
\put(-11,0){\vector(1,1){24}}
\put(13,0){\line(-1,1){10}}
\put(-1,14){\vector(-1,1){10}}

\put(64,0){
\put(-37,10){\mbox{$- \ \, A^{-1}\cdot H\Bigg($}}
\put(39,10){\mbox{$\Bigg)$}}
\put(37,0){\vector(-1,1){24}}
\put(13,0){\line(1,1){10}}
\put(27,14){\vector(1,1){10}}
}

\put(152,0){
\put(-39,10){\mbox{$= \, (q-q^{-1})\cdot H\Bigg($}}
\put(46,10){\mbox{$\Bigg)$}}
\put(32,0){\vector(0,1){24}}
\put(42,0){\vector(0,1){24}}
}
}

\end{picture}
\ee
with a requirement that the HOMFLY-PT polynomial for the unknot is equal to one. Comparing this skein relation with the Jones skein relation~\eqref{JonesSkein}, one can make sure that the HOMFLY-PT polynomial reduces to the Jones polynomial at the point $A=q^2$.
}~\cite{freyd1985new,przytycki1987kobe} which is a generalisation of the Jones polynomial. Namely, the question is whether the HOMFLY-PT polynomial detects the unknot. A. Stoimenow carried out lots of research on the Jones and HOMFLY-PT unknot problem. In 2003, he posed a theorem that there is no positive knots\footnote{Positive knots are knots having diagrams with all crossings positive.} with the trivial HOMFLY-PT polynomial~\cite{stoimenow2003positive}. Later he also proved the non-triviality of the Jones polynomial on $k$-almost positive knots\footnote{A diagram is $k$-almost positive if it has exactly $k$ negative crossings. A knot is $k$-almost positive if it has a $k$-almost positive, but no $(k-1)$-almost positive diagram.} for $k\leq 3$~\cite{stoimenow2011diagram}. A. Stoimenov also proposed a proof that there does not exist a three-strand braid whose closure gives a knot with the trivial Jones polynomial~\cite{stoimenow2017properties}. 
 
 S. Bigelow proposed a method for constructing a non-trivial knot with a trivial Jones polynomial under certain conditions on the corresponding braids~\cite{bigelow2002does}. S. Bigelow also connected the Jones problem with the problem of faithfulness of the Burau representation for
the 4-strand braid group ${\cal B}_4$ and faithfulness of the Temperley-Lieb representation, what relies the Jones problem with open classical mathematical questions. Later T. Ito proved that if the Burau representation for 4-braids is not faithful, than for any knot, there exist infinitely many mutually different knots with the same HOMFLY-PT polynomial~\cite{ito2015kernel}. However, A. Datta showed that the Burau representation of ${\cal B}_4$ is faithful almost everywhere~\cite{datta2022strong}. However some tricky cases are still open.

 Many famous mathematicians tried to solve the Jones problem. Thus, it is reasonable to apply recently discovered untested approaches for study of quantum knot invariants and to consider particular cases (for example, knots of fixed number of strands). We get use of Reshetikhin--Turaev approach~\cite{TURAEV1992865,Reshetikhin:1990pr} for construction of quantum knot invariants and different expansions of the HOMFLY polynomials -- differential~\cite{IMMM,Morozov_2016,Kononov_2016,Morozov_2018,Kameyama_2020,Morozov_2019,morozov2020kntz,BM1,arxiv.2205.12238}, character~\cite{mironov2012character} and perturbative~\cite{Alvarez1994tt, Alvarez_1997,chmutov_duzhin_mostovoy_2012, Kontsevich} expansions. We mostly analyse hypothetical 4-strand knots\footnote{In these notations, $\beta$ is a braid and its overlining means that a knot ${\cal K}:=\overline{\beta}$ is the closure of a braid $\beta$.} ${\cal K}:=\overline{\beta}$ of the unit Jones polynomials and their related knots $\tilde{\cal K}:=\overline{\beta \Delta^{2n}}$, ${\cal K}^{2k+1}=\overline{\beta^{2k+1}}$. In this study we are in particular motivated by the mentioned S. Bigelow's research. 

 This paper is organised as follows. In Section~\ref{sec:methodology}, we introduce all basic definitions and methods we use in our study of the Jones problem. We analyse the cases of two-, three- and four-strand knots in Sections~\ref{sec:two-strand},~\ref{sec:three-strand} and~\ref{sec:four-strand} correspondingly. On two and three strands, there is no non-trivial knots with the trivial Jones polynomial. An interesting case is the case of four-strand braids. We have found the HOMFLY-PT polynomial for hypothetical knots of the unit Jones polynomial. The case of a large number of strands is also discussed in Section~\ref{sec:more-strand}. In Section~\ref{sec:futher-checks}, we provide checks of the obtained HOMFLY-PT polynomials. We explore the character and perturbative expansions of these polynomials in Sections~\ref{sec:4-strandChExp},~\ref{sec:4-strandPertExp}. The character coefficients turn out to be fixed up to one coefficient $a_{[2,2]}$ which can be parameterised by an arbitrary two-component link on three strands which comes from the hypothetical four-strand braid with the trivial Jones polynomial. The verification of unitarity of the $\mathfrak{R}$-matrix representation (see Section~\ref{sec:eigenval-analysis}) turns out to be restrictive enough. Only a few two-component links can in fact parameterise the obtained HOMFLY-PT polynomial, what also gives a hint on how the hypothetical 4-strand knots of the unit Jones polynomial look like. If a knot ${\cal K}=\overline{\beta}$ of the trivial Jones polynomial does exist, then these must exist related knots $\tilde{\cal K}=\overline{\beta \Delta^{2n}}$ and ${\cal K}^{2k+1}=\overline{\beta^{2k+1}}$. In Section~\ref{sec:rel-knots}, we give some verification on the existence of such knots. Finally, we propose an effective search method for positive braids which can give the obtained hypothetical HOMFLY-PT polynomials in Section~\ref{sec:pos-braids}. 

\setcounter{equation}{0}
\section{Methodology}\label{sec:methodology}

The presentation of knots in the form of braids is an effective technique for exploring many issues in knot theory. In our work this technique will also be useful for the Jones problem. For example, in the Reshetikhin-Turaev approach due to the convenience of producing polynomial invariants.
\begin{definition}
   A braid $\beta$ on $n$ strands is a set of n non-intersecting monotonically ascending curves in $R^3$ that connect n points on the lower line to n points on the upper line. 
\end{definition}

\begin{figure}[h!]
\begin{picture}(160,150)(-400,-28)

\qbezier(-256,48)(-266,48)(-266,58)
\qbezier(-256,68)(-266,68)(-266,58)
\put(-256,68){\line(1,0){229}}

\qbezier(-17,58)(-17,68)(-27,68)
\qbezier(-17,58)(-17,48)(-27,48)

\qbezier(-256,24)(-276,24)(-276,54)
\qbezier(-256,84)(-276,84)(-276,54)
\put(-256,84){\line(1,0){229}}

\qbezier(-7,54)(-7,84)(-27,84)
\qbezier(-7,54)(-7,24)(-27,24)

\qbezier(-256,0)(-286,0)(-286,50)
\qbezier(-256,100)(-286,100)(-286,50)
\put(-256,100){\line(1,0){229}}

\qbezier(3,50)(3,100)(-27,100)
\qbezier(3,50)(3,0)(-27,0)

\qbezier(-256,-24)(-296,-24)(-296,46)
\qbezier(-256,116)(-296,116)(-296,46)
\put(-256,116){\line(1,0){229}}

\qbezier(13,46)(13,116)(-27,116)
\qbezier(13,46)(13,-24)(-27,-24)

\put(-256,48){\line(1,0){10}}
\put(-256,24){\line(1,0){10}}
\put(-256,0){\line(1,0){102}}
\put(-256,-24){\line(1,0){102}}

\put(-222,48){\line(1,0){10}}
\put(-222,24){\line(1,0){10}}
\put(-246,48){\line(1,-1){24}}
\put(-222,48){\line(-1,-1){10}}
\put(-236,34){\line(-1,-1){10}}

\put(-212,48){\line(1,-1){24}}
\put(-212,24){\line(1,1){10}}

\put(-188,48){\line(-1,-1){10}}
\put(-188,24){\line(1,0){10}}
\put(-188,48){\line(1,0){10}}
\put(-154,48){\line(-1,-1){10}}

\put(-178,48){\line(1,-1){24}}
\put(-178,24){\line(1,1){10}}

\put(-154,48){\line(1,0){20}}
\put(-154,24){\line(1,0){20}}

\put(-154,-24){\line(1,1){10}}
\put(-154,0){\line(1,-1){24}}

\put(-130,0){\line(-1,-1){10}}
\put(-130,0){\line(1,0){10}}
\put(-130,-24){\line(1,0){10}}

\put(-120,0){\line(1,-1){24}}
\put(-120,-24){\line(1,1){10}}

\put(-86,-24){\line(1,1){10}}
\put(-96,0){\line(-1,-1){10}}
\put(-96,0){\line(1,0){10}}
\put(-96,-24){\line(1,0){10}}

\put(-86,0){\line(1,-1){24}}

\put(-62,0){\line(-1,-1){10}}
\put(-134,48){\line(1,0){107}}
\put(-134,24){\line(1,0){72}}

\put(-62,0){\line(1,1){24}}
\put(-62,24){\line(1,-1){10}}
\put(-38,0){\line(-1,1){10}}

\put(-38,24){\line(1,0){11}}
\put(-38,0){\line(1,0){11}}
\put(-62,-24){\line(1,0){35}}

\end{picture}
\caption{{\small An example of four-strand braid and its closure.}}
\label{fig:braid-example}
\end{figure}
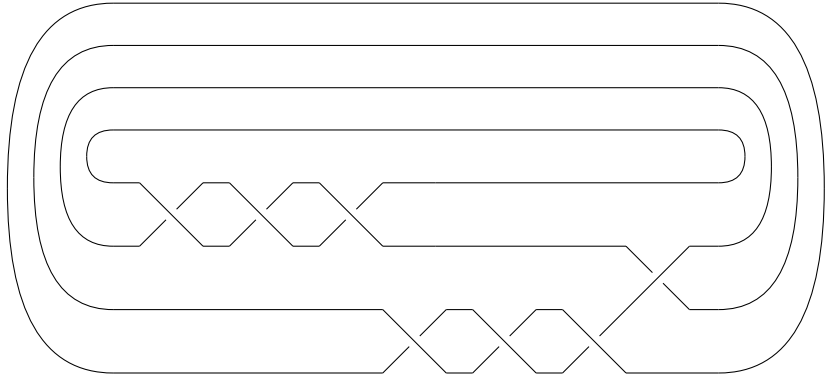

\begin{proposition}
    {\bf (Alexander theorem~{\cite{manturov2018knot}})} Every knot ${\cal{K}}$ or link ${\cal{L}}$ can be represented as a closed braid $\overline{\beta} = {\cal{K}}$.
\end{proposition}
  The closure operation is the connection of each $i$-th point on the left segment to the $i$-th point on the right segment without intertwining with other connecting curves or strands of the braid, see Fig.~\ref{fig:braid-example}. The benefit over knots is that braids have a convenient group structure. In the braid group $B_n$ the generators $\sigma_i$ correspond to the intersection of the $i$-th and $(i+1)$-th strands and satisfy the following relations.

\begin{equation}\label{braidcontext}
 \begin{cases}
   \sigma_i\sigma_{i+1}\sigma_i = \sigma_{i+1}\sigma_i\sigma_{i+1}, &   1\leq i\leq n-2\\
   \sigma_i\sigma_j = \sigma_j\sigma_i, & |i-j|>2
 \end{cases}
\end{equation}

The first expression in \eqref{braidcontext} is the famous Yang-Baxter equation. Remarkably, it corresponds to the third Reidemeister move. Besides, the second expression is called the far commutativity relation. This is also the second Reidemeister move. Thus, the relations on the braid group ensure to the invariance of the knots with respect to the Reidemeister moves.

The presentation of a knot in the braid group is matched by a product of generators $\sigma_i$ in the order of intertwining of the corresponding pairs of strands. For the four-strand braid from Fig.~\ref{fig:braid-example}, there are 3 different generators, and the representation of this braid is $\mathfrak{P}({\cal{K}}) = \sigma_1^3\sigma_3^{3}\sigma_2^{-1}$. A notation is often used where only the degrees of relevant generators in the representation of braid group are stated. $\sigma_1^{a_1}\sigma_2^{b_1}\sigma_3^{c_1}...\sigma_1^{a_n}\sigma_2^{b_n}\sigma_3^{c_n} = (a_1,b_1,c_1;...;a_n,b_n,c_n)$. For our example $\mathfrak{P}({\cal{K}}) = (3,0,3;0,-1,0)$. There are many representations of braid group. In the context of the Jones problem we are particularly interested in the Burau and the Jones (Reshetikhin-Turaev) representations. 


\subsection{The Burau representation}

The Burau representation of $\mathcal{B}_n$ has a rather obvious block-diagonal structure. The generators \eqref{burmat} differ from the unit matrix $n\times n$ on a crossing block of size $2\times2$ only in the columns whose numbers correspond to the strands to be crossed. Naturally, the generators in the Burau representation satisfies the relations of the braid group \eqref{braidcontext}.

\begin{equation}\label{burmat}\rho_i = 
\begin{pmatrix}\begin{array}{c|cc|c}
E_{i-1}  & 0 & 0 & 0 \\
\hline
0 & 1-q^2 & q^2 &  0 \\
0 &1 & 0 & 0 \\

\hline
0 & 0 & 0 & E_{n-i-1}
\end{array}
\end{pmatrix}
\end{equation}

This representation is unreduced and decomposes into one-dimensional trivial and $(n-1)$-dimensional reduced representations. The most interesting is precisely the reduced Burau representation, since it is used to determine another important --- the Alexander polynomial $Al^{\cal{K}}$. Let us give a formula for this invariant for knots on n strands\footnote{The Alexander polynomial $Al(q)$ can be defined with the use of the following skein relations: \be\label{Alexrel}
\begin{picture}(300,25)(-70,-30)

\put(0,-30){
\put(-3,10){\mbox{$ Al\Bigg( $}}
\put(40,10){\mbox{$\Bigg) $}}
\put(14,0){\vector(1,1){24}}
\put(38,0){\line(-1,1){10}}
\put(24,14){\vector(-1,1){10}}

\put(64,0){
\put(-15,10){\mbox{$- \ \, Al \Bigg($}}\put(39,10){\mbox{$\Bigg)$}}
\put(37,0){\vector(-1,1){24}}
\put(13,0){\line(1,1){10}}
\put(27,14){\vector(1,1){10}}
}

\put(152,0){
\put(-39,10){\mbox{$= \, (q-q^{-1})\cdot Al\Bigg($}}
\put(46,10){\mbox{$\Bigg)$}}
\put(32,0){\vector(0,1){24}}
\put(42,0){\vector(0,1){24}}
}
}

\end{picture}
\ee }. The values of the Alexander polynomials \eqref{Alexander} correspond to the values on "The Knot Atlas" \cite{bar2004knot}.

\begin{equation}\label{Alexander}
Al^{\cal{K}}(q) = \frac{1 - q^2}{1-q^{2n}} \, det(1-\mathfrak{B}({\cal{K}})),
\end{equation}

$\mathfrak{B}({\cal{K}})$ is reduced Burau representation for $\cal{K}$. Two-strand braids have only one generator $b = -q^2$. Generators for braids with a large braid index $n$ are defined by the following matrices of size $(n-1)\times (n-1)$.

\begin{equation}\label{b1}
b_1 =
\begin{pmatrix}
 \begin{array}{cc|c}  -q^2 & 1 & 0 \\ 0 & 1 & 0 \\ 
\hline 0 & 0 & E_{n-3} \end{array}

\end{pmatrix} \ \ \ \ \ \ \ \  b_{n-1} =
\begin{pmatrix}, 
\begin{array}{c|cc} E_{n-3} & 0 & 0 \\ \hline  0 & 1 & 0 \\
0 & q^2 & -q^2   \end{array}
\end{pmatrix}
\end{equation}

\begin{equation}\label{bi}
b_i =
\begin{pmatrix}
 \begin{array}{c|ccc|c}
E_{i-2} & 0 & 0 & 0 & 0 \\
\hline
0 & 1 & 0 & 0 & 0 \\
0 & q^2 & -q^2 & 1 & 0 \\
0 & 0 & 0 & 1 & 0 \\
\hline
0 & 0 & 0 & 0 & E_{n-i-2}
\end{array}
\end{pmatrix}, \text{ } \text{ } 2 \leq i \leq n-2
\end{equation}

\begin{definition}
  The algebraic crossing number also called the writhe number $W({\cal{K}})$ for the knot ${\cal{K}}$ is the difference of the number of positive and negative crossings on all strands in the braid representation for the knot ${\cal{K}}$. In the language of generators, it is the sum of all degrees in the braid representation.
\end{definition}

Remarkably, the Alexander polynomial does not distinguish trivial knots. For example, it is equal to one for a Conway knot $(11n34)$, as it is same for a trivial knot. This invariant is associated with the HOMFLY-PT polynomial by substitution $ A \to 1$, which will come in useful later. It is easy to see this fact by looking at the skein-relations \eqref{Alexrel} for the Alexander polynomial and comparing them to the HOMFLY-PT relations \eqref{Hrelation}.

\subsection{Reshetikhin--Turaev approach}
\label{RTapp}

\begin{definition}
    Reshetikhin--Turaev or $\mathfrak{R}$-matrix representation of the braid group ${\cal B}_n$ is a homomorphysm

    \begin{equation}
        \pi: {\cal B}_n \to {\rm End}(V_{\lambda_1} \otimes ... \otimes V_{\lambda_{n}})\,,
    \end{equation}

   \begin{equation}\label{HomoGeneratorBraidtoRMatrix}
\pi(\sigma_i) = \mathfrak{R}_i = \mathds{1} \otimes ... \otimes \mathds{1} \otimes \hat{\mathfrak{R}}_{i,i+1} \otimes \mathds{1} \otimes ... \otimes \mathds{1}\,,
\end{equation}
where $\mathds{1}$ is the identity $n \times n$ matrix, $V_{\lambda_i}$ is the finite-dimensional irreducible representation of the quantized universal enveloping algebra $U_q(\mathfrak{sl}_N)$, which is numbered by Young diagrams $\lambda_i$, $i = 1, ..., n$. $\hat{\mathfrak{R}}$ is the universal $\mathfrak{R}$-matrix for the corresponding representation of the algebra. The universal $\mathfrak{R}$-matrix has following form \cite{klimyk2012quantum}:
    \begin{equation}\label{FunR-matrix}
\hat{\mathfrak{R}}(q) = \mathbb{P}\, q^{\sum_{\alpha,\gamma} C^{-1}_{\alpha\gamma}H_{\alpha}\otimes H_{\gamma}} \prod_{\small{{\rm positive \ root} \  \beta}}\exp_q\left[(q-q^{-1})E_{\beta}\otimes F_{\beta}\right],
\end{equation}
where $\exp_q$ is the quantum exponent, $C_{\alpha\gamma}$ is the Cartan matrix, $\{E_i, F_i, H_i\}$ are generators of the algebra $U_q(sl_N)$, $\mathbb{P}$ is the permutation matrix: $\mathbb{P}(x \otimes y) = y \otimes x$.
\end{definition}

Defining the representation of the braid group, $\mathfrak{R}$-matrices satisfy the Yang--Baxter equation:
\begin{equation}\label{Yang-Baxter-equation}
\mathfrak{R}_{i}\mathfrak{R}_{i+1}\mathfrak{R}_i = \mathfrak{R}_{i+1}\mathfrak{R}_i\mathfrak{R}_{i+1}\,.
\end{equation}

\begin{definition}

In the Reshetikhin--Turaev approach, the colored HOMFLY-PT polynomial is defined as the quantum trace from the Reshetikhin--Turaev representation of a braid $\beta$, the closure of which gives a knot $\cal{K}$ or a link $\cal{L}$:
\begin{equation}\label{HOMFLY-PTformula}
H^{\cal{K}}_R(q,A=q^N) = {\rm Tr}_q\left(\prod_s \mathfrak{R}_{\mu(s)}^{\{R-T\}}\right)\,.
\end{equation}
In this formula, $s$ enumerates the interweaves in the braid representation. And $\mu(s)$ fixes the generator (the strands that are intertwined). Young diagram $R$ corresponds to a certain representation $U_q(\mathfrak{sl}_N)$, ${\rm Tr}_q$ denotes the quantum trace.  
\end{definition}

We are interested in the fundamental representation $(\lambda_i = \square)$ of the quantum algebra $U_q(\mathfrak{sl}_N)$. In this case, the Reshetikhin-Turaev representation is reduced to the Jones representation, which defines the fundamental HOMFLY-PT polynomial and, consequently, the fundamental Jones ($N = 2$) polynomial of interest \cite{jones1987hecke}.

For the fundamental representation of $U_q(\mathfrak{sl}_2)$ the universal $\mathfrak{R}$-matrix has following form:
\begin{equation}\label{FunR-matrixN2}
\hat{\mathfrak{R}}(q) = \mathbb{P}\, q^{\frac{1}{2} H\otimes H} \left(\mathds{1} + (q - q^{-1})E\otimes F\right)\,,
\end{equation}
where $\mathbb{P},\, E,\, F,\, H$ are taken in the fundamental representation of quantum algebra $U_q(\mathfrak{sl}_2)$:
\begin{equation}\label{P-matrix}
\mathbb{P} =
\begin{pmatrix}
 \begin{array}{cccc}  1&0 & 0&0   \\  0&0 & 1&0\\0& 1 & 0&0 \\ 0 & 0 & 0& 1      
 \end{array}
 \end{pmatrix};
\end{equation}

\begin{equation}\label{Sl2}
E =
\begin{pmatrix}
 \begin{array}{cc}  0 & 1   \\  0& 0 
 \end{array}
 \end{pmatrix}; \text{ } \text{ } \text{ } \text{ } \text{ } \text{ }
 F =
\begin{pmatrix}
 \begin{array}{cc}  0 & 0   \\  1& 0 
 \end{array}
 \end{pmatrix}; \text{ } \text{ } \text{ } \text{ } \text{ } \text{ }
 H =
\begin{pmatrix}
 \begin{array}{cc}  1 & 0   \\  0 & -1 
 \end{array}
 \end{pmatrix}.
\end{equation}
The quantum trace for the fundamental representation of $U_q(\mathfrak{sl}_2)$ has the following form:
\begin{equation}\label{QTr}
 {\rm Tr}_q\left(f\right) = {\rm Tr}\left(\hat{\mathfrak{K}}^{\otimes n}f\right), \ \ \ \ \ \ \ \hat{\mathfrak{K}} = {\rm diag}(q,  q^{-1})\,.
\end{equation}

\subsection{Connection of Jones and Burau representations for $\mathcal{B}_4$ group}

\begin{definition}
    A representation $T: G_n \to GL(N)$ is called \textit{faithful} if all elements $\gamma$ of the group $G_n$ are represented by different images of the linear mapping $T$. Accordingly, a representation T is \textit{unfaithful} if there exists such different elements $\beta_1, \beta_2$ of group ${\cal B}_n$ that $T(\beta_1) = T(\beta_2)$.
\end{definition}

In the early 2000s, Stephen Bigelow in \cite{bigelow2002does} associated the Jones problem ,the faithfulness question for the Jones representation (the fundamental Reshetikhin-Turaev representation) and the faithfulness problem for the Burau representation for the braid group on four strands ${\cal B}_4$.

\begin{proposition}\label{Big1}
    {\bf (S. Bigelow \cite{bigelow2002does})} If the Jones representation of the braid group is unfaithful then there exists a way to construct a nontrivial knot with a trivial HOMFLY-PT polynomial. Accordingly, the Jones polynomial for such a knot will also be trivial.
\end{proposition}

\begin{proposition}\label{Big2}
    {\bf (S. Bigelow \cite{bigelow2002does})} The Burau representation for ${\cal B}_4$ and the Jones representation for the braid group on four strands are both faithful or unfaithful.
\end{proposition}

That said, the question of the faithfulness of the Burau representation is currently open only for ${\cal B}_4$. For braid groups on an other number strands it is known that the Burau representation of ${\cal B}_n$ are faithful for $n \le 3$ and unfaithful for $n \ge 5$. The faithfulness problem for the Reshetikhin-Turaev representation and, consequently, the question of the faithfulness of the Jones representation are also open.

To illustrate proposition \ref{Big2} consider the decomposition of the Jones representation of ${\cal B}_4$ into a direct sum of irreducible representation $V_{\lambda}$:
\begin{equation}\label{JoneshowBurau}
Jones_4 = V_{[4]} \oplus V_{[3,1]} \oplus V_{[2,2]} \oplus  V_{[2,1,1]} \oplus  V_{[1,1,1,1]} \,.
\end{equation}
All terms in \eqref{JoneshowBurau} are Burau matrices. The representations $V_{[4]}, V_{[1,1,1,1]}$ are one-dimensional and correspond to Burau matrices for two strands. More interesting are the three-dimensional representations $V_{[3,1]}$, $V_{[2,1,1,1]}$ which corresponds to Burau matrices for four strands and the two-dimensional representation $V_{[2,2]}$ which correspond Burau matrices for three strands. The complete $\mathfrak{R}$-matrix \eqref{HomoGeneratorBraidtoRMatrix} in the basis of irreducible representations has block-diagonal structure. So the following $\mathfrak{R}$-matrices are from \cite{mironov2012character} which are actually diagonal blocks in the full matrix. Thus all of them are equivalent to reducible Burau matrices of the corresponding dimensions.
\begin{equation}\label{R-matrix1}
\mathfrak{R}_1^{[2,2]} = \mathfrak{R}_3^{[2,2]} =
\begin{pmatrix}
 \begin{array}{cc}  q & 0   \\  0& -\frac{1}{q} 
 \end{array}

\end{pmatrix}; \text{ } \text{ } \text{ } \text{ } \text{ } \text{ }
\mathfrak{R}_2^{[2,2]} =
\begin{pmatrix}
 \begin{array}{cc}  -\frac{1}{q^2[2]_q} & -\frac{\sqrt{[3]_q}}{[2]_q}  \\   -\frac{\sqrt{[3]_q}}{[2]_q}& \frac{q^2}{[2]_q} 
 \end{array}

\end{pmatrix};
\end{equation}

\begin{equation}\label{R-matrix2}
\mathfrak{R}_1^{[3,1]} =
\begin{pmatrix}
 \begin{array}{ccc}  q & 0 & 0 \\ 0 & q & 0 \\0 & 0& -\frac{1}{q} 
 \end{array}

\end{pmatrix}; \text{ } \text{ } \text{ } 
\mathfrak{R}_2^{[3,1]} =
\begin{pmatrix}
 \begin{array}{ccc}  q & 0 & 0 \\ 0 & -\frac{1}{q^2[2]_q} & -\frac{\sqrt{[3]_q}}{[2]_q} \\0&-\frac{\sqrt{[3]_q}}{[2]_q} &  \frac{q^2}{[2]_q}
 \end{array}

\end{pmatrix};
\end{equation}
\begin{equation}\label{R-matrix3}
\mathfrak{R}_3^{[3,1]} =
\begin{pmatrix}
 \begin{array}{ccc}  -\frac{1}{q^3[3]_q} &  -\frac{[2]_q\sqrt{q^2 + q^{-2}}}{[3]_q} & 0 \\ -\frac{[2]_q\sqrt{q^2 + q^{-2}}}{[3]_q} & \frac{q^3}{[3]_q} & 0 \\0 & 0& q 
 \end{array}

\end{pmatrix}.
\end{equation}

The matrices for the $[2,1,1]$ representation are obtained from the matrices for the $[3,1]$ representation by substituting $q \to -\frac{1}{q}$. It is not difficult to see that such a substitution corresponds to a negative inverse matrix for $[3,1]$, what can also be rewritten in terms of the inversion of the parameter $q$:
\begin{equation}\label{R-matrix4}
   \mathfrak{R}_i^{[2,1,1]}(q) = \mathfrak{R}_i^{[3,1]}\left(-\frac{1}{q}\right) = - \left(\mathfrak{R}_i^{[3,1]}\right)^{-1}(q)= -\mathfrak{R}_i^{[3,1]}(q^{-1})\,.
\end{equation}
The Burau representation for braid group on three strings and the Jones representation for $B_3$ are also isomorphic:
\begin{equation}\label{JoneshowBurau3}
Jones_3 = V_{[3]} \oplus V_{[2,1]} \oplus V_{[1,1,1]}\,. 
\end{equation}
The representations $V_{[3]}$ and $V_{[1,1,1]}$ are one-dimensional. And the representation $V_{[2,1]}$ is two-dimensional and corresponds to Burau matrices for two strands. The $\mathfrak{R}-$matrices for this representation are identical to the matrices from \eqref{R-matrix1}.

\begin{equation}\label{R-matrix[2,1]}
    \mathfrak{R}_1^{[2,1]} = \mathfrak{R}_1^{[2,2]}; \ \ \ \ \ \ \ \mathfrak{R}_2^{[2,1]} = \mathfrak{R}_2^{[2,2]}\,.
\end{equation}
It is easy to see that such $\mathfrak{R}$-matrices are equivalent to Burau matrices. After a proper rescaling of Burau matrices~\eqref{b1},~\eqref{bi} of the corresponding size and inversion of $q$, one can find a basis transformation which gives the $\mathfrak{R}$-matrices~\eqref{R-matrix1}--\eqref{R-matrix3}.




Interesting corollary of the equivalence of the fundamental Reshetikhin-Turaev and Burau representations for braid groups on 3 and 4 strands is the possibility to define the Alexander polynomial in terms of the $\mathfrak{R}$-matrices \eqref{AlexanderhowR3} and \eqref{AlexanderhowR4}: 
\begin{equation}\label{AlexanderhowR3}
Al^{\cal{K}}(q) = q^{2 + W}\,\frac{1 - q^2}{1-q^6} \, \det(E - q^{-W}\mathfrak{R}(K)) \quad \text{for 3-braids},
\end{equation}

\begin{equation}\label{AlexanderhowR4}
Al^{\cal{K}}(q) = q^{3 + W}\,\frac{1 - q^2}{1-q^8} \, \det(E - q^{-W}\mathfrak{R}(K))  \quad \text{for 4-braids}.
\end{equation}
The HOMFLY-PT polynomial for knots on three or four strands can be similarly determined via the Burau matrices.

We remark that we do not know the inverse theorem to the theorem \ref{Big1}, so the faithfulness of the Burau representation for ${\cal B}_4$ does not mean that the Jones problem is closed.

\subsection{Character expansion}

It is much more convenient to consider not the complete $\mathfrak{R}$-matrices, but their diagonal blocks \eqref{R-matrix1} and \eqref{R-matrix2}, since in this case the problem is reduced to the analysis of representations of lower dimension. Then the HOMFLY-PT polynomials can be character decomposed in some basis $S_{\lambda}$. This approach is known as \textit{character expansion}.
\begin{equation}\label{characterexpansion}
H^{\cal{K}} = \sum_{\lambda} h^{\lambda} S^{*}_{\lambda}\,.
\end{equation}
The index $\lambda \in \square^{\otimes n}$ is responsible for decomposition into irreducible representations, which are numbered by Young diagrams, $h^{\lambda}$ is the trace function of the reduced Reshetikhin-Turaev representation of the braid $\beta \in {\cal B}_{n}$.

\begin{definition}\label{defShur}
    
$S_{\lambda}$ is the Schur polynomial at the special point $p_k = \frac{A^k - A^{-k}}{q^k - q^{-k}}$, which is given by the hooks formula \eqref{kruk}:
 \begin{equation}\label{kruk}
S_Q = \prod_{(i,j)\in Q}\frac{\{Aq^{i-j}\}}{\{q^{c_{i,j}}\}} = \prod_{(i,j)\in Q}\frac{[N+i-j]_q}{[c_{i,j}]_q} \begin{picture}(105,30)(-35,-25)
\put(0,0){\line(1,0){70}}
\put(0,-10){\line(1,0){70}}
\put(0,-20){\line(1,0){60}}
\put(0,-30){\line(1,0){40}}
\put(0,-40){\line(1,0){20}}

\put(0,0){\line(0,-1){40}}
\put(10,0){\line(0,-1){40}}
\put(20,0){\line(0,-1){40}}
\put(30,0){\line(0,-1){30}}
\put(40,0){\line(0,-1){30}}
\put(50,0){\line(0,-1){20}}
\put(60,0){\line(0,-1){20}}
\put(70,0){\line(0,-1){10}}
\put(15,-15){\makebox(0,0)[cc]{\textbf{$\star$}}}
\put(15,-25){\makebox(0,0)[cc]{\textbf{$\star$}}}
\put(15,-35){\makebox(0,0)[cc]{\textbf{$\star$}}}

\put(25,-15){\makebox(0,0)[cc]{\textbf{$\star$}}}
\put(35,-15){\makebox(0,0)[cc]{\textbf{$\star$}}}
\put(45,-15){\makebox(0,0)[cc]{\textbf{$\star$}}}
\put(55,-15){\makebox(0,0)[cc]{\textbf{$\star$}}}
\put(15,5){\makebox(0,0)[cc]{$i$}}
\put(55,5){\makebox(0,0)[cc]{$k$}}
\put(-5,-15){\makebox(0,0)[cc]{$j$}}
\put(-5,-35){\makebox(0,0)[cc]{$l$}}
\put(25,-37){\makebox(0,0)[lc]{$c_{i,j}=k-i+l-j+1$}}
\end{picture}
\end{equation}
where $\{x\} = x - x^{-1}, [x]_q = \frac{q^x - q^{-x}}{q - q^{-1}}, A = q^N,  c_{i,j}$ is a hook length (number of cubes in the fragment of a Young diagram $Q$ in the form of a hook). 
\end{definition}
On the presented example in picture \eqref{kruk}, the hook for the cell with coordinates $i = 2$ and $j = 2$ is marked by stars $\star$. It is easy to see that for the presented Young diagram $c_{i,j} = 7$.

Polynomial invariants account for the normalisation of Schur polynomials by $S_{[1]} = \frac{A-A-A^{-1}}{q-q^{-1}}$. It is necessary to normalise the polynomial of the trivial knot by one. That is the Schur functions from \eqref{characterexpansion} are expressed through the Schur functions from \eqref{kruk} as follows:
\begin{equation}\label{NormShur}
S^{*}_Q = S_Q\cdot  \frac{q-q^{-1}}{A-A^{-1}}\,.
\end{equation}

For knots on 2, 3 and 4 strands respectively, the following formulas for the HOMFLY-PT polynomials in the fundamental representation hold:

\begin{equation}\label{HOMFLYhowR2}
H^{\cal{K}}(q, A) = A^{-W}\left(q^W S^{*}_{[2]} + (-q)^{-W}S^{*}_{[1,1]}\right)\,,
\end{equation}

\begin{equation}\label{HOMFLYhowR3}
H^{\cal{K}}(q, A) = A^{-W}\left(q^W S^{*}_{[3]} + \tr(\mathfrak{R}^{[2,1]})S^{*}_{[2,1]} + (-q)^{-W}S^{*}_{[1,1,1]}\right)\,,
\end{equation}

\begin{equation}\label{HOMFLYhowR4}
H^{\cal{K}}(q, A) = A^{-W}\left(q^W S^{*}_{[4]} + \tr(\mathfrak{R}^{[3,1]})S^{*}_{[3,1]}+ \tr(\mathfrak{R}^{[2,2]})S^{*}_{[2,2]}+ \tr(\mathfrak{R}^{[2,1,1]})S^{*}_{[2,1,1]} + (-q)^{-W}S^{*}_{[1,1,1,1]}\right)\,.
\end{equation}
The Jones polynomial is obtained by substituting $A \to q^2$. In this case, the formulas \eqref{HOMFLYhowR2} - \eqref{HOMFLYhowR4} are greatly simplified, since with this substitution the Schur polynomials are nullified in diagrams with more than two rows. For two-row diagrams, there exist isomorphism relationships for the corresponding representations:
\begin{equation}
   [2,2] = [1,1] = \varnothing; \ \ \ \ \ [3,1] = [2]; \ \ \ \ \ [2,1] = [1]\,.
\end{equation}
In the formalism of Schur diagrams, the isomorphism relationships take the following form:
\begin{equation}
   S^{*}_{[2,2]} = S^{*}_{[1,1]} = \frac{q}{1+q^2}\,; \ \ \ \ \ S^{*}_{[3,1]} = S^{*}_{[2]} = \frac{(1 - q + q^2)(1 + q + q^2)}{q+q^3}\,; \ \ \ \ \ S^{*}_{[2,1]} = S^{*}_{[1]} = 1\,.
\end{equation}

\subsection{Differential expansion}

The interesting behaviour of the HOMFLY-PT polynomials when $A \to q$ is substituted. It is easy to see from the formula for the Schur functions \eqref{kruk} that only the first Schur functions take a value different from zero $S^{*}_{[n]} = 1$. The coefficient in front of them is also equal to one. Then $H^{\cal{K}} (A = q) = 1$. Therefore, we can rewrite our invariant through another polynomial $C^{\cal{K}}$ so that the reduction to unity at $A \to q$ is explicitly visible in this representation. This technique is called \textit{differential expansion} \cite{arxiv.2205.12238}.

\begin{equation}\label{difexp}
H^{\cal{K}} = 1 + (A^{2} - q^{2})\cdot C^{\cal{K}}\,.
\end{equation}
The quadraticity in the multiplier is due to the fact that the HOMFLY-PT polynomial in our notations is a function on $A^2, q^2$. We note that such conditions of triviality of invariants under certain substitutions of polynomial variables can strengthen the differential expansion. In addition there are HOMFLY-PT symmetries that can increase the number of multipliers at $C^{\cal{K}}$. Thus there exists a symmetry related to the transpose of Young diagrams $R$:

\begin{equation}\label{symmetryH}
H^{\cal{K}}_R(q,A) = H^{\cal{K}}_{R^T}( q^{-1},A)\,.
\end{equation}

In our work, we explore an exclusively fundamental representation. Obviously transposition of it will also be the fundamental representation $(\square^T = \square)$. Then this symmetry takes the following form:
\begin{equation}\label{symmetryHfun}
H^{\cal{K}}(q, A) = H^{\cal{K}}(q^{-1}, A)\,.
\end{equation}
Now it is easy to see that when we substitute $A \to q^{-1}$, the invariants stop distinguishing between knots because they become trivial. So, the differential expansion can be transformed:
\begin{equation}\label{difexp2}
H^{\cal{K}} = 1 + (A^{2} - q^{2})(A^{2} - q^{-2})C^{\cal{K}}\,.
\end{equation}

\subsection{Perturbative expansion}\label{sec:PertExpDef}

To research the properties of the HOMFLY-PT polynomial $H^{\cal{K}}_R(A, q)$, it is often useful to expand it into a formal series on the variable $\hbar$ by substituting $q \to e^{\hbar/2}, A \to e^{N \hbar/2}$. This approach is called \textit{perturbative expansion} or \textit{Vassiliev invariant expansion}. A wonderful property of the perturbative expansion is to separate the dependence of the polynomial on the knot and the representation.

\begin{definition}
    Vassiliev invariants are the coefficients ${\cal V}_{n,m}^{\cal{K}}$ in the perturbative expansion of the HOMFLY-PT polynomial for a knot $\cal{K}$;

    \begin{equation}\label{VasInvStruc}
H_R^{\cal{K}}(q, A)|_{q \to e^{\hbar/2}, A \to e^{N\hbar/2}}  = \sum_{n=0}^{\infty}\left(\sum_{m=1}^{{\rm dim}(\mathbb{G}_n)}{\cal V}_{n,m}^{\cal{K}}{\cal G}^R_{n,m}\right)\hbar^n,
\end{equation}
where $R$ indicates the representation of the algebra $\mathfrak{sl}_N$,  ${\cal G}^R_{n,m}$ are called group factors, ${\rm dim}(\mathbb{G}_n)$ is the number of independent group factors at order $n$.
\end{definition}

Group factors are known up to 13th order \cite{lanina2021chern}, \cite{lanina2022implications}. The first few terms of the perturbative expansion are as follows:
\begin{equation}\label{ExampleVasInvStruc}
H_R^{\cal K}(q, A) = 1 + \hbar^2{\cal V}_{2,1}^{\cal{K}}{\cal G}^R_{2,1} + \hbar^3{\cal V}_{3,1}^{\cal{K}}{\cal G}^R_{3,1} + \hbar^4({\cal V}_{4,1}^{\cal{K}}{\cal G}^R_{4,1} + {\cal V}_{4,2}^{\cal{K}}{\cal G}^R_{4,2} + {\cal V}_{4,3}^{\cal{K}}{\cal G}^R_{4,3}) + O(\hbar^5)\,.
\end{equation}
The first invariant is zero for all knots what can be easily seen from the general form of the differential expansion for the HOMFLY-PT polynomial \eqref{difexp2}. For the unknot all Vassiliev invariants are zero. 
Note that in the perturbative expansion of the HOMFLY polynomial~\eqref{VasInvStruc}, not all group factors are independent. Some of them factor into products of and group factors of lower levels. Vassiliev invariants inherit the same property. For example, at the 4-th order
\begin{equation}
    {\cal G}_{4,1}={\cal G}^2_{2,1}\quad \text{and} \quad {\cal V}_{4,1}={\cal V}^2_{2,1}\,.
\end{equation}
This sort of dependence is exponential. I.e. all independent group factors and Vassiliev invariants (called {\it primitive}) are included into $\log (H^{\cal K}_R)$.

In what follows, we consentrate only on the case of fundamental representation $R$. At the first glance, one can think that the condition that a polynomial in $A^2$, $q^2$ possesses the perturbative expansion~\eqref{VasInvStruc} is rather restrictive. However, it can be shown that this condition is equivalent to the condition that a polynomial in $A^2$, $q^2$ can be represented in the form of differential expansion
\begin{equation}
    H^{\cal{K}}(q, A) = 1 + (A^{2} - q^{2})(A^{2} - q^{-2})\cdot C^{\cal{K}}(q, A) \quad \text{with} \quad C^{\cal{K}}(q, A)=C^{\cal{K}}(q^{-1}, A)\,.
\end{equation}

\setcounter{equation}{0}
\section{Analysis of knots of fixed braid width}\label{sec:fixed-braid-width-analysis}
In this section we analyse possible knots of fixed braid width that can give the unit Jones polynomial. First, using Reshetikhin--Turaev approach and character expansion, one can easily observe that there is no non-trivial 2-strand knot with the unit Jones polynomial, see Section~\ref{sec:two-strand}. Second, in Section~\ref{sec:three-strand}, we provide a proof that there is no non-trivial 3-strand knot with $J^{\cal K}=1$ and $H^{\cal K} \neq 1$. Using analogous methods, in Section~\ref{sec:four-strand}, we prove that there is no non-trivial 4-strand knot with $Al^{\cal K}=J^{\cal K}=1$ and $H^{\cal K} \neq 1$, and further study the case of $J^{\cal K}=1$, $Al^{\cal K} \neq 1$, $H^{\cal K} \neq 1$. In Section~\ref{sec:more-strand}, we provide a similar analysis for braids of larger number of strands.

\subsection{The case of 2 strands}\label{sec:two-strand}
At first, let us analyse the simplest case of two-strand knots.

\begin{proposition}
    There does not exist a non-trivial two-strand knot with the trivial Jones polynomial.   
\end{proposition}
$\Delta$ The Jones polynomial for two-strand knots is given by
\begin{equation}
    J^{\cal K}(q)={\rm Tr}_q\, \left(\mathfrak{R}^{2n+1}\right)= q^{-4n-2}\left(q^{2n+1}S^*_{[2]}+(-q)^{-2n-1}S^*_{[1,1]}\right)\,.
\end{equation}
 Then, the Jones polynomial is
\begin{equation}
    J^{\cal K}(q)=\frac{q^{-2n +1} + q^{-2n - 1} + q^{-2n -3} -q^{-6n - 3}}{q+q^{-1}}\,,
\end{equation}
and it is easy to see that it is equal to one if and only if $n=0$ or $n=-1$ what corresponds to two different diagrams of the unknot. All other two-strand knots have non-trivial Jones polynomials. $\Box$

\subsection{The case of 3 strands}\label{sec:three-strand}
Reshetikhin--Turaev approach does not look so simple for three-strand knots. In this case, there are two different $\mathfrak{R}$-matrices $\mathds{1}\otimes 
\mathfrak{R}$ and $\mathfrak{R}\otimes \mathds{1}$ which do not commute, and thus, cannot be diagonalised in the same basis. Thus, in this section we utilise the differential expansion for the HOMFLY-PT polynomial. Another useful tool for our analysis is the Morton-Franks-Williams inequality relating the $A$-span of the HOMFLY-PT polynomial ${\rm span}_A H^{\cal K}(A,q)$ with the braid index $\mathfrak{b}^{\cal K}$ of a knot ${\cal K}$.
\begin{definition}
    The $A$-span of the HOMFLY-PT polynomial is the difference between maximum and minimum degrees in $A$.
\end{definition}
\begin{definition}
    The braid index $\mathfrak{b}^{\cal K}$ of a knot ${\cal K}$ is the least number of strings needed to make a closed braid representation of a knot ${\cal K}$.
\end{definition}
\begin{proposition}
    {\bf (Morton~\cite{morton1986seifert,morton1988polynomials}, Franks and Williams~\cite{franks1987braids})} \\ Let ${\rm span}_A H^{\cal K}(q, A)$ be the $A$-span of the HOMFLY-PT polynomial of a knot ${\cal K}$, and $\mathfrak{b}^{\cal K}$ be the braid index. Then the Morton-Franks-Williams inequality holds:
    \begin{equation}\label{MFW}
\frac{1}{2}\left({\rm span}_A H^{\cal K}(q, A)\right)  \leq \mathfrak{b}^{\cal K} - 1
\end{equation}
\end{proposition}
Now, we are ready to state the following theorem.
\begin{proposition}
\label{3strand}
    There does not exist a non-trivial three-strand knot with the trivial Jones polynomial and the non-trivial HOMFLY-PT polynomial. 
\end{proposition}
$\Delta$ Assume that there exists a non-trivial knot with $J^{\cal K}=1$ and $H^{\cal K}\neq 1$. The substitution $A=q^2$ reduces the HOMFLY-PT polynomial to the Jones polynomial. Thus, the differential expansion~\eqref{difexp2} gets a stronger form
\begin{equation}\label{difexp3br}
    H^{\cal K}(q, A) = 1 + (A^2 - q^2)(A^2 - q^{-2})(A^2 - q^4)(A^2 - q^{-4})\cdot C^{\cal K}\,.
    \end{equation}
For knots of the braid index equal to three, the Morton-Franks-Williams inequality~\eqref{MFW} implies that
\begin{equation}
    \frac{1}{2}{\rm span}_A H^{\cal K}(q, A)  \leq 2\,.
\end{equation}
In that way, the HOMFLY-PT polynomial of our hypothetical knots must have the $A$-span equal to 4 or less. Cyclotomic function $C^{\cal K}$ is a polynomial in $A^2$, $q^2$. Due to the differential expansion~\eqref{difexp3br}, the $A$-span of the HOMFLY-PT polynomial is not less than 6 in general case. By fixing the function $C^{\cal K}$, we can reduce only one power of $A^2$ in the HOMFLY-PT polynomial. To decrease the $A$-span, the highest or the lowest power of $A^2$ must be vanished. In other words, choosing a proper cyclotomic function $C^{\cal K}$, one can decrease the $A$-span of the HOMFLY polynomial~\eqref{difexp3br} only to 6. Thus, we arrive to the contradiction, and such knots do not exist. $\Box$  
\\ \\
The proved theorem serves as a preliminary to the 4-strand case. Actually, A. Stoimenow stated the general theorem.
\begin{proposition}
    {\bf (Stoimenow \cite{stoimenow2017properties})} There does not exist a non-trivial 3-strand knot with the trivial Jones polynomial.
\end{proposition}

\subsection{The case of 4 strands}\label{sec:four-strand}
Finally, let us study currently open question on existence of a 4-strand knot with the unit Jones polynomial. 
\begin{theorem}
    There does not exist a non-trivial 4-strand knot with the trivial Jones and Alexander polynomials and the non-trivial HOMFLY-PT polynomial.
\end{theorem}
$\Delta$ A proof is analogous to the one of Proposition \ref{3strand}. Now, the differential expansion becomes stronger as we demand the substitution $A=1$ to reduce the HOMFLY-PT polynomial to the unit Alexander polynomial:
\begin{equation}
H^{\cal K} = 1 + (A^2 - 1)(A^2 - q^2)(A^2 - q^{-2})(A^2 - q^4)(A^2 - q^{-4})\cdot C^{\cal K}\,.
\end{equation}
Then, on the one hand, the minimal number of different $A^2$ powers in the HOMFLY-PT polynomial is equal to $5$. On the other hand, the Morton-Franks-Williams inequality~\eqref{MFW} implies that the number of different $A^2$ powers must be not greater than $4$. So, we again arrive to the contradiction, and we conclude that 4-strand knots with $J^{\cal K}=Al^{\cal K}=1$ and $H^{\cal K}\neq 1$ do not exist. $\Box$
\\ \\
Now, we consider hypothetical 4-strand knots with the trivial Jones polynomial and the non-trivial Alexander and HOMFLY-PT polynomials. 

\begin{theorem}
\label{HOM4}
    If there exist 4-strand knots with the trivial Jones polynomial and the non-trivial HOMFLY-PT polynomial, then the HOMFLY-PT polynomials of such knots have the following differential expansion:
    \begin{equation}\label{difF}
H^{\cal K}(q, A) = 1 + (A^2 - q^2)(A^2 - q^4)(A^2 - q^{-2})(A^2 - q^{-4})\cdot F^{\cal K}(q, A)
\end{equation}
    with the following cyclotomic functions:
    \begin{equation}\label{CycF}
    \begin{aligned}
        F^{{\cal K}_m} &= -\sum^m_{j=0}\frac{[2j+4]_q[j+3]_q[j+1]_q}{[3]_q[4]_q}A^{-2j - 8}\,,\\
        \bar{F}^{\bar{{\cal K}}_m} &= -\sum^m_{j=0}\frac{[2j+4]_q[j+3]_q[j+1]_q}{[3]_q[4]_q}A^{2j}\,, 
    \end{aligned}
    \end{equation}
with an arbitrary non-negative integer number $m$.
\end{theorem}
$\Delta$ As we have already noted in the proof of Proposition \ref{3strand}, the polynomial~\eqref{difF} has the $A$-span equal to 6 or more in general case. The Morton-Franks-Williams inequality~\eqref{MFW} demands that this value must be at most $6$. Thus, our goal is to decrease the $A$-span by a proper choice of $F^{\cal K}$-polynomial. Let us use the following notations for short:
\begin{equation}\label{ShDE}
H^{\cal K} = 1 +\left[ C_0 + C_2 A^2 + C_4 A^4 + C_6A^6 + C_8 A^8\right]F^{\cal K}
\end{equation}
with
\begin{equation}
\begin{aligned}
    C_0 &= C_8 = 1\,, \\
    C_2 &= C_6 = - (q^4+q^2+q^{-2}+q^{-4})\,, \\
    C_4 &= q^6+q^2+q^{-2}+q^{-6}\,.
\end{aligned}
\end{equation}
One can easily guess that $A^0$- or $A^8$-monomials in square brackets can be vanished by the first term of~\eqref{ShDE}, which is just unity, choosing
\begin{equation}
    \bar{F}^{\bar{{\cal K}}_0} = -1 \quad \text{or} \quad F^{{\cal K}_0} = -A^{-8}\,.
\end{equation}
Actually, there is a whole series of such polynomials~\eqref{CycF}. One can check that in these polynomials all odd $A^2$ powers are mutually absorbed. Let us explain this fact for $\bar{F}^{\bar{{\cal K}}_m}$. Denote
\begin{equation}\label{Fpolynomial}
\bar{F}^{\bar{{\cal K}}_m} =\sum^m_{j=0}\phi_{2j} A^{2j}\,.
\end{equation}
The Morton-Franks-Williams condition can be straightforwardly checked. The first several equations on $\phi_{2n}$ are
\begin{equation}\label{eq1}
\begin{aligned}
    C_0\phi_0 + 1 = 0\,,\\
    C_2\phi_0 + C_0\phi_2 = 0\,,\\
    C_4\phi_0 + C_2\phi_2 + C_0\phi_4 = 0\,,\\
    C_6\phi_0 + C_4\phi_2 + C_2\phi_4 + C_0\phi_6 = 0\,.
\end{aligned}
\end{equation}
For all other $\phi_{2j}$-coefficients the constraint on $A$-span of the HOMFLY-PT polynomial~\eqref{difF} is fulfiled due to the following relation:
\begin{equation}\label{eq2}
C_8\phi_{2i} + C_6\phi_{2i+2} + C_4\phi_{2i+4} + C_2\phi_{2i+6} + C_0\phi_{2i+8} = 0
\end{equation}
that can be checked for all $\phi_{2j}$ in a general form. Note that the coefficients of equations~\eqref{eq1},~\eqref{eq2} form the lower triangular matrix with the determinant equal to $C_0^{m+1}=1$. Thus, the corresponding linear system of equations is non-degenerate, and its solution is unique.  

Another series of $F^{\cal{K}}$-polynomials has the form
\begin{equation}
F^{{\cal K}_m} =\sum^m_{j=0}\phi_{2j}
A^{-2j-8}\,
\end{equation}
with the same coefficients $\phi_{2j}$. Thus, the constraint that the $A$-span of the HOMFLY-PT polynomial~\eqref{difF} is at most $6$ also holds. $\Box$
\begin{remark}
    If polynomials $F^{{\cal K}_m}$ correspond to some knots ${\cal K}_m$, then polynomials $\bar{F}^{\bar{{\cal K}}_m}$ correspond to knots $\bar{{\cal K}}_m$ mirror to ${\cal K}_m$.
\end{remark}
$\Delta$ Indeed, if
\begin{equation}
H^{{\cal K}_m}(q, A) = 1 + (A^2 - q^2)(A^2 - q^4)(A^2 - q^{-2})(A^2 - q^{-4})\cdot F^{{\cal K}_m}(q, A)\,,
\end{equation}
then
\begin{equation}
  H^{\bar{{\cal K}}_m}(q, A)=H^{{\cal K}_m}(q^{-1}, A^{-1})=1+(A^2 - q^2)(A^2 - q^4)(A^2 - q^{-2})(A^2 - q^{-4})\cdot \overbrace{A^{-8}F^{{\cal K}_m}(A^{-1},q^{-1})}^{=\bar{F}^{\bar{{\cal K}}_m}(q, A)}\,.\quad \Box
\end{equation}
If knots with $J^{\cal K}(q)=1$ exist, then $J^{\bar{{\cal K}}}(q)=J^{\cal K}(q^{-1})=1$, and mirror knots $\bar{{\cal K}}$ also possess the trivial Jones polynomials. It is a good sign that the obtained HOMFLY-PT polynomials correspond to knots and their mirror images. Also note that $H^{{\cal K}_m}(q, A)$ and $H^{\bar{{\cal K}}_m}(q, A)$ are different functions, so that the hypothetical {\it knots ${\cal K}_m$ must be chiral} \footnote{A chiral knot is a knot that is not equivalent to its mirror image.}. In what follows, we focus only on the family ${\cal K}_m$. 
\begin{remark}
    Polynomials $F^{{\cal K}_m}$ possess the invariance under inversion of $q$:
    \begin{equation}
        F^{{\cal K}_m}(q, A)=F^{{\cal K}_m}(q^{-1}, A)\,,
    \end{equation}
    what corresponds to the symmetry of the fundamental HOMFLY-PT polynomials
    \begin{equation}
        H^{\cal K}(q, A)=H^{\cal K}(q^{-1}, A)\,.
    \end{equation}
\end{remark}

\subsection{The case of a large number of strands}\label{sec:more-strand}

The case with a large number of strands in a braid is more complex and includes many representations that cannot be reduced to Burau matrices. Analysing these representations for faithfulness is a challenging task (see attempts to do so for a 5-braid here \cite{kasahara2008remarks}). Nevertheless, we can consider the limiting case where the MFW inequality \eqref{MFW} must be sharp, by analogy with Theorem \ref{HOM4}. 

To do so, we need to impose conditions on the HOMFLY-PT polynomial as the number of strands in the braid increases. As additional conditions, it is reasonable to consider trivialisation of Reshetikhin--Turaev invariants for a fixed $\mathfrak{sl}_n$ algebra, as discussed in Section \ref{RTapp}. These invariants are obtained from the HOMFLY-PT polynomial by specialising the variable $A$ to $q^n$. The Jones polynomial is therefore the $\mathfrak{sl}_2$ Reshetikhin--Turaev invariant. With a larger number of strands, in addition to trivialising the Jones polynomial, it may also be possible to trivialise $\mathfrak{sl}_3$ invariants and so on. In the case of odd number of strands we add one additional condition of the triviality of the Alexander polynomial.



\begin{theorem}[\textbf{case of even $\mathfrak{b}^{\cal K}$}]
Suppose there exist $\mathfrak{b}^{\cal K}$-strand knots with $\mathfrak{sl}_n$ Reshetikhin-Turaev invariants being trivial up to $n \le \frac{\mathfrak{b}^{\cal K}}{2}$. Then the HOMFLY-PT polynomials of such knots have the following differential expansion:
    \begin{equation}\label{difF}
H^{\cal K}(q, A) = 1 +  F^{\cal K} (q, A) \cdot\prod^{\mathfrak{b}^{\cal K}}_{\substack{ k = - \mathfrak{b}^{\cal K} \\  k \neq 0}}\left(A^2 - q^{2k}\right)
\end{equation}
where cyclotomic functions are:
    \begin{equation}\label{CycF}
    \begin{aligned}
        F^{{\cal K}_m} &= -\sum^m_{j=0}\frac{[2j+\mathfrak{b}^{\cal K}]_q \, [\frac{\mathfrak{b}^{\cal K}}{2}]_q \, [j,{\mathfrak{b}^{\cal K}}]_q}{[\mathfrak{b}^{\cal K}]!_q \, [j + \frac{\mathfrak{b}^{\cal K}}{2}]_q} \, A^{-2j - 2 \mathfrak{b}^{\cal K}}\,,
        \\
        \bar{F}^{\bar{{\cal K}}_m} &=  -\sum^m_{j=0}
        \frac{[2j+\mathfrak{b}^{\cal K}]_q \, [\frac{\mathfrak{b}^{\cal K}}{2}]_q \, [j,{\mathfrak{b}^{\cal K}}]_q}{[\mathfrak{b}^{\cal K}]!_q \, [j + \frac{\mathfrak{b}^{\cal K}}{2}]_q} \, A^{2j}\,, 
    \end{aligned}
    \end{equation}
with q-factorial $[j]!_q := [j]_q \cdot [j - 1]_q \cdot ... \cdot [2]_q \cdot [1]_q$, rising q-factorial $[j,k]_q := [j]_q \cdot [j + 1]_q \cdot ... \cdot [j + k -1]_q$ and $m$ is an arbitrary non-negative integer number.
\end{theorem}

\begin{theorem}[\textbf{case of odd $\mathfrak{b}^{\cal K}$}]
Suppose there exist $\mathfrak{b}^{\cal K}$-strand knots with $\mathfrak{sl}_n$ Reshetikhin-Turaev invariants being trivial up to $n \le \frac{\mathfrak{b}^{\cal K}-1}{2}$ and trivial Alexander polynomial $Al^{\cal K}(q) = 1$. Then the HOMFLY-PT polynomials of such knots have the following differential expansion:
\begin{equation}\label{difF}
H^{\cal K}(q, A) = 1 +  F^{\cal K} (q, A) \cdot\prod^{\mathfrak{b}^{\cal K} - 1}_{k = - \mathfrak{b}^{\cal K} + 1}\left(A^2 - q^{2k}\right)
\end{equation}
    with the following cyclotomic functions:
    \begin{equation}\label{CycF}
    \begin{aligned}
        F^{{\cal K}_m} &= -\sum^m_{j=0}\frac{[j,\mathfrak{b}^{\cal K}{-}1]_q}{[\mathfrak{b}^{\cal K}{-}1]_q!} \, A^{-2j - 2\mathfrak{b}^{\cal K}}\,,\\
        \bar{F}^{\bar{{\cal K}}_m} &=  -\sum^m_{j=0}\frac{[j,\mathfrak{b}^{\cal K}{-}1]_q}{[\mathfrak{b}^{\cal K}{-}1]_q!} \, A^{2j}\,. 
    \end{aligned}
    \end{equation}
\end{theorem}

The proofs of these two theorems are completely analogous to the proof of Theorem \ref{HOM4}.

\setcounter{equation}{0}
\section{Verifications of obtained 4-strand HOMFLY-PT polynomials}\label{sec:futher-checks}
In the previous section we have obtained the HOMFLY-PT polynomials for hypothetical 4-strand knots with the unit Jones polynomial~\eqref{difF},~\eqref{CycF}. A question we raise and give some answers to in this section is how one can get sure that the obtained polynomials are indeed the HOMFLY-PT polynomials of some knots. We analyse the character expansion and perturbative structure of the obtained hypothetical HOMFLY polynomials in Sections~\ref{sec:4-strandChExp},~\ref{sec:4-strandPertExp}. In Section~\ref{sec:eigenval-analysis}, we check the unitarity property of the $\mathfrak{R}$-matrix representation of the hypothetical 4-strand knots at certain values of the parameter $q$.
This check prohibits many 2-links that could parameterise the obtained polynomials~\eqref{difF},~\eqref{CycF} from the previous section. In Section~\ref{sec:rel-knots}, we investigate related knots $\tilde{\cal K}:=\overline{\beta \Delta^{2n}}$, ${\cal K}^{2k+1}=\overline{\beta^{2k+1}}$ for new conditions on the hypothetical knot $\cal{K}$ with the unit Jones polynomial. In Section~\ref{sec:pos-braids}, we propose a convenient method for computerised search of a 4-strand knot with the trivial Jones polynomial and the non-trivial HOMFLY polynomial.

\subsection{Character expansion}\label{sec:4-strandChExp}
If the found polynomials $H^{{\cal K}_m}$ really correspond to some knots ${\cal K}_m$, then the 4-strand character expansion~\eqref{HOMFLYhowR4} for these HOMFLY polynomials must hold:
\begin{equation}\label{Hforcoef}
H^{\cal K}(q, A) = A^{-W}\left(a_{[4]} S^{*}_{[4]} + a_{[3,1]}S^{*}_{[3,1]}+ a_{[2,2]}S^{*}_{[2,2]}+ a_{[2,1,1]}S^{*}_{[2,1,1]} + a_{[1,1,1,1]}S^{*}_{[1,1,1,1]}\right)\,.
\end{equation}
We compare the differential expansion of $H^{{\cal K}_m}$~\eqref{difF} with its character expansion~\eqref{Hforcoef}. Polynomials $H^{{\cal K}_m}$ given by~\eqref{difF} and the first formula in~\eqref{CycF} have four $A^2$ powers, what gives a system of 4 equations on 5 $a_Y$-coefficients of character expansion~\eqref{Hforcoef}. Note that from the 4-strand character expansion~\eqref{HOMFLYhowR4}, we know that $a_{[4]}=q^W$ and $a_{[1,1,1,1]}=-q^{-W}$, and comparison with the differential expansion~\eqref{difF} with cyclotomic functions $F^{{\cal K}_m}$~\eqref{CycF} fixes the writhe number $W=2m+5$. As we have 4 constraints on 5 coefficients, we set $a_{[2,2]}(q)$ to be a free parameter. As a result, we get
\begin{equation}\label{a-coef}
\begin{aligned}
    a_{[4]} &= q^{2m + 5}, \\ a_{[1,1,1,1]} &= - q^{-2m - 5}\,, \\
    a_{[2,1,1]} &= -\frac{1 + q^2 \, - \, q^{2 + 2m}\left(1+q^2+q^4+q^6+q^8\right) \, + \, a_{[2,2]}q^{11 + 4m}}{q^{9 + 4m}(1-q+q^2)(1+q+q^2)}\,, \\
    a_{[3,1]} &= -\frac{a_{[2,2]}q^2 \, + \,  q^{2m+3}\left(1+q^2+q^4+q^6+q^8\right) \, - \, q^{11 + 4m}\left(1+q^2\right) }{(1-q+q^2)(1+q+q^2)}\,.
\end{aligned}
\end{equation}
A simple observation is the fact that our hypothetical 4-strand knots of the trivial Jones polynomial have odd writhe number, and the minimal writhe number is $W_{\rm min}=5$.

We express $a_{[2,1,1]}$ and $a_{[3,1]}$ through $a_{[2,2]}$ because $a_{[2,2]}$ is the trace of $2\times 2$ matrix. It is easier to analyse that $a_{[2,1,1]}$ and $a_{[3,1]}$ corresponding to some $3\times 3$ matrices. Moreover, the coefficient $a_{[2,2]}$ for a knot ${\cal K}=\overline{\beta}$ on 4 strands is equal to the coefficient $a_{[2,1]}$ for a knot ${\cal K}'=\overline{\beta'}=\overline{\beta\big|_{\sigma_3\rightarrow \sigma_1}}$ on 3 strands. I.e. in 4-strand braid $\beta$ we can substitute all $\sigma_3$ generators to $\sigma_1$ generators what brings us to 3-strand braid $\beta'$; for an example see Fig.~\ref{fig:four-three-strands}. This follows from the fact that $\mathfrak{R}_1^{[2,2]} = \mathfrak{R}_3^{[2,2]} = \mathfrak{R}_1^{[2,1]}$ and $\mathfrak{R}_2^{[2,2]} = \mathfrak{R}_2^{[2,1]}$, see~\eqref{R-matrix1},~\eqref{R-matrix[2,1]}. In other words, considering coefficient $a_{[2,2]}$, we effectively deal with 3-strand braids. The closure of these 3-strand braids actually give a 2-component link due to the fact that the writhe number must be odd. 

Thus, the coefficients $a_{[2,2]}$ are parameterised by 2-component links on 3-strands with the writhe number $W=2m+5$, $m=0,1,2,\dots$ The simplest example is provided by the torus link $T[2,6]$ that can be represented as the closure of 3-strand braid $(6,-1)$, see Fig.~\ref{fig:four-three-strands}. This diagram of the torus link $T[2,6]$ has $W=5$, $m=0$. Then, the coefficients of the character expansion of the HOMFLY-PT polynomial for the hypothetical 4-strand knot are
\begin{equation}
\begin{aligned}
    a_{[2,2]} &= a_{[2,1]} = \tr\left(\left(\mathfrak{R}^{[2,1]}_1\right)^6\left(\mathfrak{R}^{[2,1]}_2\right)^{-1}\right) = \frac{(1-q^2)(1+q^4)(1+q^8)}{q^7}\,, \\
    a_{[4]} &= q^{5}, \quad a_{[3,1]} = \frac{-a_{[2,2]} - q -q^3 - q^5 -q^7 + q^{11}}{(1-q+q^2)(1+q+q^2)}q^2\,, \\
    a_{[1,1,1,1]} &= -q^{-5}, \quad a_{[2,1,1]} = -\frac{q^{11}a_{[2,2]} - q^{10} -q^8 - q^6 -q^4 + 1}{q^9(1-q+q^2)(1+q+q^2)}\,,
\end{aligned}
\end{equation}
and the resulting HOMFLY-PT polynomial is
\begin{equation}\label{HP}
    \small{H^{\cal K} = \frac{1}{A^8q^6}(A^6q^2-A^4+A^2q^2+A^2q^4-A^4q^4+A^6q^4-q^6 - 2A^4q^6 +A^2q^8-A^4q^8+A^6q^8+A^2q^{10} + A^6q^{10} -A^4q^{12})}\,.
\end{equation}
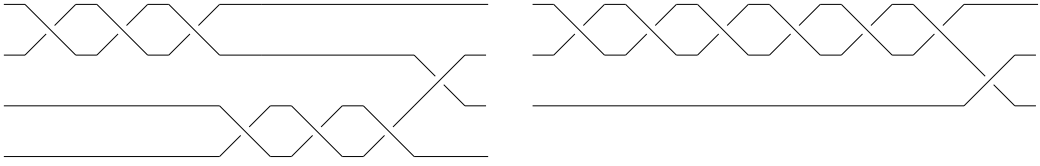
\begin{figure}[h!]
\scalebox{0.8}{
\begin{picture}(160,90)(-320,-35)
\put(-6,48){\line(1,0){10}}
\put(-6,24){\line(1,0){10}}
\put(-6,0){\line(1,0){204}}

\put(28,48){\line(1,0){10}}
\put(28,24){\line(1,0){10}}
\put(4,48){\line(1,-1){24}}
\put(28,48){\line(-1,-1){10}}
\put(14,34){\line(-1,-1){10}}

\put(38,48){\line(1,-1){24}}
\put(38,24){\line(1,1){10}}

\put(62,48){\line(-1,-1){10}}

\put(62,24){\line(1,0){10}}
\put(72,48){\line(1,-1){24}}
\put(72,24){\line(1,1){10}}
\put(96,48){\line(-1,-1){10}}

\put(62,48){\line(1,0){10}}
\put(96,48){\line(1,0){10}}
\put(96,24){\line(1,0){10}}
\put(106,24){\line(1,1){10}}
\put(106,48){\line(1,-1){24}}
%
\put(130,48){\line(-1,-1){10}}
\put(130,48){\line(1,0){10}}
\put(130,24){\line(1,0){10}}
%
\put(140,48){\line(1,-1){24}}
\put(140,24){\line(1,1){10}}

\put(174,24){\line(1,1){10}}
\put(164,48){\line(-1,-1){10}}
\put(164,48){\line(1,0){10}}
\put(164,24){\line(1,0){10}}

\put(174,48){\line(1,-1){24}}

\put(198,48){\line(-1,-1){10}}
\put(198,48){\line(1,0){35}}

\put(198,0){\line(1,1){24}}
\put(198,24){\line(1,-1){10}}
\put(222,0){\line(-1,1){10}}

\put(222,24){\line(1,0){10}}
\put(222,0){\line(1,0){10}}


\put(-256,48){\line(1,0){10}}
\put(-256,24){\line(1,0){10}}
\put(-256,0){\line(1,0){102}}
\put(-256,-24){\line(1,0){102}}

\put(-222,48){\line(1,0){10}}
\put(-222,24){\line(1,0){10}}
\put(-246,48){\line(1,-1){24}}
\put(-222,48){\line(-1,-1){10}}
\put(-236,34){\line(-1,-1){10}}

\put(-212,48){\line(1,-1){24}}
\put(-212,24){\line(1,1){10}}

\put(-188,48){\line(-1,-1){10}}
\put(-188,24){\line(1,0){10}}
\put(-188,48){\line(1,0){10}}
\put(-154,48){\line(-1,-1){10}}

\put(-178,48){\line(1,-1){24}}
\put(-178,24){\line(1,1){10}}

\put(-154,48){\line(1,0){20}}
\put(-154,24){\line(1,0){20}}

\put(-154,-24){\line(1,1){10}}
\put(-154,0){\line(1,-1){24}}

\put(-130,0){\line(-1,-1){10}}
\put(-130,0){\line(1,0){10}}
\put(-130,-24){\line(1,0){10}}

\put(-120,0){\line(1,-1){24}}
\put(-120,-24){\line(1,1){10}}

\put(-86,-24){\line(1,1){10}}
\put(-96,0){\line(-1,-1){10}}
\put(-96,0){\line(1,0){10}}
\put(-96,-24){\line(1,0){10}}

\put(-86,0){\line(1,-1){24}}

\put(-62,0){\line(-1,-1){10}}
\put(-134,48){\line(1,0){107}}
\put(-134,24){\line(1,0){72}}

\put(-62,0){\line(1,1){24}}
\put(-62,24){\line(1,-1){10}}
\put(-38,0){\line(-1,1){10}}

\put(-38,24){\line(1,0){10}}
\put(-38,0){\line(1,0){10}}
\put(-62,-24){\line(1,0){35}}

\end{picture}
}
\caption{{\small An example of four-strand braid $\beta$ is on the left hand side. The three-strand braid $\beta'$ on the right hand side is obtained from the braid $\beta$ by replacement of intersections between the third and the fourth strands to the same number of intersections between the first and the second strands. After closure, this braid $\beta'$ gives two-component torus link $T[2,6]$.}}
\label{fig:four-three-strands}
\end{figure}

\subsubsection*{Properties of coefficients of character expansion}
Let us find out properties of $a_Y$-coefficients of the character expansion~\eqref{Hforcoef}. We get use of these properties in Section~\ref{sec:rel-knots}. 

Let $\mathcal{R}$ be the $\mathfrak{R}$-matrix representation in fundamental representation of our hypothetical four-strand braid $\beta$ with odd writhe number, and $\mathcal{R}^Y$, $Y=[1,1,1,1]$, $[2,1,1]$, $[3,1]$, $[2,2]$, $[4]$, are diagonal blocks of $\mathcal{R}$ in the basis of irreducible representations $Y$. The trace of any matrix is invariant under transposition. In particular, $\tr\left(\left({\cal R}^Y\right)^T\right)=\tr \left({\cal R}^Y\right)$. One can also check that $\mathfrak{R}$-matrices $\mathfrak{R}_i^Y$~\eqref{R-matrix1},~\eqref{R-matrix2},~\eqref{R-matrix3},~\eqref{R-matrix4} obey the relation $\mathfrak{R}_i^Y(q^{-1})=\left(\mathfrak{R}_i^Y(q)\right)^{-1}$. These two facts imply that
\begin{equation}\label{R-inverse}
    \tr \left({\cal R}(q^{-1})\right)=\tr\left( {\cal R}^{-1}(q)\right)\,.
\end{equation}
For any $2\times 2$ matrix $M$, $\tr\left( M\right)=\det \left(M\right)\cdot\tr \left(M^{-1}\right)$ (in particular, this property holds for $\mathfrak{R}_i^{[2,2]}$). For any $i=1,2,3$, $\det\left( \mathfrak{R}_i^{[2,2]}\right)=-1$, and the writhe number of the braid $\beta$ is odd, so that for $a_{[2,2]}=\tr\left( {\cal R}^{[2,2]}\right)$ the following equality holds:
\begin{equation}\label{a22-rel}
    a_{[2,2]}(q^{-1})=-a_{[2,2]}(q)\,.
\end{equation}
Due to the connection between $\mathfrak{R}_i^{[2,1,1]}$ and $\mathfrak{R}_i^{[3,1]}$~\eqref{R-matrix4}, we get the relation between $a_{[3,1]}$ and $a_{[2,1,1]}$:
\begin{equation}
    a_{[3,1]}(q^{-1})=-a_{[2,1,1]}(q)\,,
\end{equation}
what can also be seen from explicit expressions for these coefficients~\eqref{a-coef}. In total, we get the following general relation below. This fact is obvious and commonly known.
\begin{proposition}
For the coefficients of character expansion of the 4-strand HOMFLY polynomial~\eqref{Hforcoef}, the following relation holds:
 \begin{equation}\label{a-rel}
    a_Y(q^{-1})=-a_{Y^T}(q)\,.
\end{equation}   
\end{proposition}

Note that relation~\eqref{a-rel} serves as a verification that~\eqref{Hforcoef} with odd writhe number $W$ has the symmetry $H^{\cal K}(q, A)=H^{\cal K}(q^{-1}, A)$. It can be also shown that the similar property holds for the eigenvalues $\lambda_i^Y$ of ${\cal R}^Y$.
\begin{proposition}
Given $\mathfrak{R}$-matrix representation of a 4-strand braid, its eigenvalues satisfy the following relation: 
    \begin{equation}\label{LambdaRel}
    \forall\, i\;\;\; \exists\, j\,:\; \lambda_i^Y(q)=-\lambda_j^{Y^T}(q^{-1})\,,
\end{equation}
where a Young diagram $Y$ enumerates blocks of this $\mathfrak{R}$-matrix representation corresponding to irreducible representations.
\end{proposition}
$\Delta$ Let us prove it for $Y=[2,1,1]$. For this purpose, we construct the system of equations on $\left( \lambda_1^{[2,1,1]}, \, \lambda_2^{[2,1,1]}, \, \lambda_3^{[2,1,1]} \right)\,$, which determines them. This system can be obtained as follows. First, for $\mathfrak{R}$-matrix representation ${\cal R}$ for a knot ${\cal K}$
\begin{equation}
    \tr\left( {\cal R}^{[2,1,1]}\right)=a_{[2,1,1]}=\lambda_1^{[2,1,1]} + \lambda_2^{[2,1,1]} + \lambda_3^{[2,1,1]}\,.
\end{equation}
Second, as we have derived~\eqref{R-inverse}
\begin{equation}
    \tr \left(\left({\cal R}^{[2,1,1]}(q)\right)^{-1}\right)=\tr \left({\cal R}^{[2,1,1]}(q^{-1})\right)=a_{[2,1,1]}(q^{-1})=\left(\lambda_1^{[2,1,1]}\right)^{-1} + \left(\lambda_2^{[2,1,1]}\right)^{-1} + \left(\lambda_3^{[2,1,1]}\right)^{-1}\,.
\end{equation}
Third, for the determinant we get
\begin{equation}
    \det \left({\cal R}^{[2,1,1]}\right)=\left(\det\left( \mathfrak{R}^{[2,1,1]}_1\right)\right)^W=q^{-W}= \lambda_1^{[2,1,1]} \cdot \lambda_2^{[2,1,1]} \cdot \lambda_3^{[2,1,1]}
\end{equation}
due to the fact that $\det \left(\mathfrak{R}^{[2,1,1]}_1\right)=\det\left( \mathfrak{R}^{[2,1,1]}_2\right)=\det \left(\mathfrak{R}^{[2,1,1]}_3\right)=q^{-1}$.

Thus, the system of equations on $\lambda_i^{[2,1,1]}(q)$ is
\begin{equation}\label{lambda4}
 \begin{cases}
   \lambda_1^{[2,1,1]}(q) + \lambda_2^{[2,1,1]}(q) + \lambda_3^{[2,1,1]}(q) = a_{[2,1,1]}(q) \phantom{\Bigg |}\\
   \lambda_1^{[2,1,1]}(q)\cdot \lambda_2^{[2,1,1]}(q)\cdot \lambda_3^{[2,1,1]}(q) = q^{-W} \phantom{\Bigg |}\\
   \left(\lambda_1^{[2,1,1]}(q)\right)^{-1} + \left(\lambda_2^{[2,1,1]}(q)\right)^{-1} + \left(\lambda_3^{[2,1,1]}(q)\right)^{-1} = a_{[2,1,1]}(q^{-1}) \phantom{\Bigg |}
 \end{cases}
\end{equation}
and on $\lambda_i^{[3,1]}(q^{-1})$ is
\begin{equation}\label{lambda3}
 \begin{cases}
   \lambda_1^{[3,1]}(q^{-1}) + \lambda_2^{[3,1]}(q^{-1}) + \lambda_3^{[3,1]}(q^{-1}) = a_{[3,1]}(q^{-1})=-a_{[2,1,1]}(q) \phantom{\Bigg |}\\
   \lambda_1^{[3,1]}(q^{-1})\cdot \lambda_2^{[3,1]}(q^{-1})\cdot \lambda_3^{[3,1]}(q^{-1}) = -q^{-W} \phantom{\Bigg |}\\
   \left(\lambda_1^{[3,1]}(q^{-1})\right)^{-1} + \left(\lambda_2^{[3,1]}(q^{-1})\right)^{-1} + \left(\lambda_3^{[3,1]}(q^{-1})\right)^{-1} = a_{[3,1]}(q)=-a_{[2,1,1]}(q^{-1}) \phantom{\Bigg |}
 \end{cases}
\end{equation}
Comparing~\eqref{lambda4} and~\eqref{lambda3}, we conclude that sets of solutions of these systems coincide, and arrive to the particular case of~\eqref{LambdaRel}. Note that system~\eqref{lambda4} is equivalent to the characteristic equation 
\begin{equation}\label{LambdaCharEq}
\lambda^3-\lambda^2\cdot a_{[2,1,1]}(q)+\lambda\cdot q^{-W} a_{[2,1,1]}(q^{-1})-q^{-W}=0\,,
\end{equation}
and thus, system~\eqref{lambda4} has only one solution. For the eigenvalues $\lambda_i^{[2,2]}$ the statement~\eqref{LambdaRel} is obvious due to~\eqref{a22-rel} and the system
\begin{equation}\label{lambda22-sys}
    \begin{cases}
        \lambda_1^{[2,2]}(q)+\lambda_2^{[2,2]}(q)=a_{[2,2]}(q) \\
        \lambda_1^{[2,2]}(q)\cdot \lambda_2^{[2,2]}(q)=-1
    \end{cases}
\end{equation}
Its solution is also unique due to its equivalence to the quadratic characteristic equation.  $\Box$

\subsection{Perturbative expansion}\label{sec:4-strandPertExp}
As it has been already stated in Section~\ref{sec:PertExpDef}, the constraint that a polynomial $H^{\cal K}(q, A)$ in $A^2$, $q^2$ has correct perturbative expansion is equivalent to the constraint that this polynomial $H^{\cal K}(q, A)$ can be represented in the form of differential expansion~\eqref{difexp2} with polynomial cyclotomic functions and that it satisfies the relation $H^{\cal K}(q, A)=H^{\cal K}(q^{-1}, A)$. This is actually the case of our polynomials $H^{{\cal K}_m}$~\eqref{difF},~\eqref{CycF}. However, let us show that they indeed can be perturbative expanded and analyse the perturbative structure of these hypothetical HOMFLY-PT polynomials obtained in Section~\ref{sec:four-strand}. First, note that due to the fact that
\begin{equation}
    H^{{\cal K}_m}-1\sim (A^2-q^2)(A^2-q^4)(A^2-q^{-2})(A^2-q^{-4})
\end{equation}
the perturbative expansion of $H^{{\cal K}_m}$ is of the form
\begin{equation}
    H^{{\cal K}_m}=1+O(h^4)\,.
\end{equation}
I.e. Vassiliev invariants of orders $2$ and $3$ vanish. Second, it is convenient to analyse primitive Vassiliev invariants in the recently discovered basis of the HOMFLY-PT group factors ${\cal C}_{[k]}$ known up to the 13th order~\cite{lanina2021chern,lanina2022implications}:
\begin{equation}\label{PrimVasFind}
\begin{aligned}  &\log(H_R^{\mathcal{K}})=1+\hbar^2\left(v^{\mathcal{K}}_{[2],0}\right)_2\mathcal{C}_{[2]}^{R}+\hbar^3\left(v^{\mathcal{K}}_{[2],1}\right)_3N\mathcal{C}_{[2]}^{R}+\hbar^4\left(\left(v^{\mathcal{K}}_{[2],2}\right)_4N^2\mathcal{C}_{[2]}^{R}+\left(v^{\mathcal{K}}_{[4],0}\right)_4\mathcal{C}_{[4]}^{R}\right)+\\
&+\hbar^5\left(\left(v^{\mathcal{K}}_{[2],1}\right)_5 N\mathcal{C}_{[2]}^{R}+\left(v^{\mathcal{K}}_{[2],3}\right)_5 N^3\mathcal{C}_{[2]}^{R}+\left(v^{\mathcal{K}}_{[4],1}\right)_5 N\mathcal{C}_{[4]}^{R}\right)+\\
&+\hbar^6\left(\left(v^{\mathcal{K}}_{[2],2}\right)_6 N^2\mathcal{C}_{[2]}^{R}+\left(v^{\mathcal{K}}_{[2],4}\right)_6 N^4\mathcal{C}_{[2]}^{R}+\left(v^{\mathcal{K}}_{[4],0}\right)_6 \left(\mathcal{C}_{[4]}^{R}+\frac{1}{6}\left(\mathcal{C}_{[2]}^R\right)^2\right)+\left(v^{\mathcal{K}}_{[4],2}\right)_6 N^2\mathcal{C}_{[4]}^{R}+\left(v^{\mathcal{K}}_{[6],0}\right)_6\mathcal{C}_{[6]}^{R}\right)+O(\hbar^7),
\end{aligned}
\end{equation}
We have found the primitive Vassiliev invariants for knots ${\cal K}_m$ up to 13th order by the comparison of the perturbative expansion~\eqref{PrimVasFind} for representation $R=\Box$ with the series expansion of logarithm of~\eqref{Hforcoef} with coefficients~\eqref{a-coef}. We list here only the Vassiliev invariants of up to 6th order:
\begin{equation}
\begin{aligned}
    \left(v_{[2],0}^{\cal K}\right)_2 &=0\,, \quad \left(v_{[2],1}^{\cal K}\right)_3 =0\,, \\
    \left(v_{[2],2}^{\cal K}\right)_4 &= -\frac{5}{12}(m+1) (m+2) (m+3) (m+4)\,,\quad \left(v_{[4],0}^{\cal K}\right)_4 = 4 (m+1) (m+2) (m+3) (m+4)\,, \\
    \left(v_{[2],3}^{\cal K}\right)_5 &= \frac{1}{6} \left(6 \left(v_{[2],1}^{\cal K}\right)_5+2 m^5+25 m^4+120 m^3+275 m^2+298 m+120\right)\,,\\ 
    \left(v_{[4],1}^{\cal K}\right)_5 &= -\frac{4}{5} \left(15 \left(v_{[2],1}^{\cal K}\right)_5+4 m^5+50 m^4+240 m^3+550 m^2+596 m+240\right), \\
    \left(v_{[4],0}^{\cal K}\right)_6 &= \frac{1}{24} \left(-144 \left(v_{[2],2}^{\cal K}\right)_6+144 \left(v_{[2],4}^{\cal K}\right)_6+18 m^6+270 m^5+1655 m^4+5300 m^3+9307 m^2+8410 m+3000\right)\,, \\ 
    \left(v_{[4],2}^{\cal K}\right)_6 &= \frac{1}{6} \left(-72 \left(v_{[2],4}^{\cal K}\right)_6-6 m^6-90 m^5-533 m^4-1580 m^3-2449 m^2-1870 m-552\right)\,, \\
    \left(v_{[6],0}^{\cal K}\right)_6 &= 10 (m+1) (m+2) (m+3) (m+4) \left(2 m^2+10 m+5\right)\,.
\end{aligned}
\end{equation}
The Vassiliev invariants of order $n$ are polynomials in $m$ of order $n$. Note that some of the Vassiliev invariants are not fixed because information on other colored HOMFLY polynomials is needed.

\subsection{Analysis of eigenvalues}\label{sec:eigenval-analysis}
If the polynomials~\eqref{Hforcoef} with the coefficients of the character expansion given by~\eqref{a-coef} are the HOMFLY-PT polynomials of some knots ${\cal K}_m$, then these character expansion coefficients must come from $\mathfrak{R}$-matrix representation of the corresponding braids $\beta_m$. Three-strand~\eqref{R-matrix[2,1]},~\eqref{R-matrix1} and four-strand~\eqref{R-matrix2},~\eqref{R-matrix3},~\eqref{R-matrix4} $\mathfrak{R}$-matrices possess unitarity at the special values of $q$.
\begin{theorem}
    {\bf \cite{wenzl1992unitary,kolganov2023large}} 3-strand $\mathfrak{R}$-matrices in fundamental representation~\eqref{R-matrix[2,1]},~\eqref{R-matrix1} are unitary at the points $q=\exp\frac{2\pi i}{k}$ for $k \geq 6, \; \text{or}\; k \leq -6, \; \text{or}\; \frac{6}{1+3 z}\leq k\leq \frac{6}{-1+3 z}\; \forall\; z\in \mathbb{Z}/\{0\}$. 4-strand $\mathfrak{R}$-matrices in fundamental representation~\eqref{R-matrix2},~\eqref{R-matrix3},~\eqref{R-matrix4} are unitary at the points $q=\exp\frac{2\pi i}{k}$ for $k \geq 8, \; \text{or}\; k \leq -8, \; \text{or}\; \frac{8}{1+4 z}\leq k\leq \frac{8}{-1+4 z}\; \forall\; z\in \mathbb{Z}/\{0\}$.
\end{theorem}
$\Delta$ Give a proof for three-strand $\mathfrak{R}_i^{[2,1]}$. 
Obviously, the diagonal matrix $\mathfrak{R}_1^{[2,1]}$ is unitary at all values $q=\exp\frac{2\pi i}{k}$ for real $k$. The $\mathfrak{R}$-matrix $\mathfrak{R}_2^{[2,1]}$ is obtained from $\mathfrak{R}_1^{[2,1]}$ via the change of basis~\cite{mironov2012character}:
\begin{equation}\label{UR}
\mathfrak{R}_2^{[2,1]} = U
\begin{pmatrix}
 \begin{array}{cc}  q & 0   \\  0& -\frac{1}{q} 
 \end{array}
\end{pmatrix}U^T\,; \ \ \ \ \ U =\begin{pmatrix}
 \begin{array}{cc}  C & S   \\ -S & C 
 \end{array}
\end{pmatrix};
\ \ \ \ \ C = \frac{1}{[2]_q}\,, \ \ \ S = \frac{\sqrt{[3]_q}}{[2]_q}\,.
\end{equation}
Then, $\mathfrak{R}_2^{[2,1]}$ is unitary if the matrix $U$ is unitary:
\begin{equation}\label{UR2}
 UU^{\dagger} =
\begin{pmatrix}
 \begin{array}{cc}  \frac{|q|^2(1 + |1 + q^2 + q^{-2}|)}{(1 + q^2)(1 + (q^*)^2)} &  \frac{2i|q|^2{\rm Im}\left(\sqrt{1 + q^2 + q^{-2}}\right)}{(1 + q^2)(1 + (q^*)^2)}  \\ -\frac{2i|q|^2 {\rm Im}\left(\sqrt{1 + q^2 + q^{-2}}\right)}{(1 + q^2)(1 + (q^*)^2)} &  \frac{|q|^2(1 +|1 + q^2 + q^{-2}|)}{(1 + q^2)(1 + (q^*)^2)} 
 \end{array}
\end{pmatrix}\,,
\end{equation}
where we denote complex conjugation by $*$. Unitarity of $U$ means that $UU^\dagger$ must be the unit matrix. At $q=\exp\frac{2 \pi i}{k}$, it is easy to see that the diagonal elements of $UU^\dagger$ are units. The non-diagonal elements equal to zero if
\begin{equation}
    1 + q^2 + q^{-2} \ge 0\,.
\end{equation}
At the point $q=\exp\frac{2 \pi i}{k}$ we get the following condition and the allowed real values of $k$:
\begin{equation}
    \cos \frac{4\pi}{k} \ge -\frac{1}{2} \quad \Rightarrow \quad k \geq 6, \quad \text{or}\quad k \leq -6, \quad \text{or}\quad \frac{6}{1+3 z}\leq k\leq \frac{6}{-1+3 z}\quad \forall\; z\in \mathbb{Z}/\{0\}\,.
\end{equation}
For 4-strand $\mathfrak{R}$-matrices the condition of unitarity is stronger:
\begin{equation}\label{4-str-restr}
    \cos \frac{4\pi}{k} \ge 0 \quad \Rightarrow \quad k \geq 8, \quad \text{or}\quad k \leq -8, \quad \text{or}\quad \frac{8}{1+4 z}\leq k\leq \frac{8}{-1+4 z}\quad \forall\; z\in \mathbb{Z}/\{0\}\,. \quad \Box
\end{equation}

\begin{corollary}
    Since the $\mathfrak{R}$-matrices \eqref{R-matrix1}--\eqref{R-matrix3} are unitary at the particular values $q=\exp\frac{2\pi i}{k}$ for $k \geq 8, \; \text{or}\; k \leq -8, \; \text{or}\; \frac{8}{1+4 z}\leq k\leq \frac{8}{-1+4 z}\; \forall\; z\in \mathbb{Z}/\{0\}$, the Jones representation for any braid is also unitary at these points. Hence, its eigenvalues lie on the unit circle, $|\lambda_i| = 1$.  
\end{corollary} 
Although we cannot find out the particular form of the matrix representation of the hypothetical braids from the obtained polynomials~\eqref{difF},~\eqref{CycF}, we can find the corresponding eigenvalues. We have already written down the system of equations on the eigenvalues $\lambda_i^{[2,1,1]}$, see~\eqref{lambda4}. For representation $[3,1]$, the system takes the similar form:
\begin{equation}\label{lambdaeq}
 \begin{cases}
   \lambda_1(q) + \lambda_2(q) + \lambda_3(q) = a_{[3,1]}(q)\\
   \lambda_1(q) \cdot \lambda_2(q) \cdot \lambda_3(q) = (-q)^W\\
   \lambda_1^{-1}(q) + \lambda_2^{-1}(q) + \lambda_3^{-1}(q) = a_{[3,1]}(q^{-1})
 \end{cases}
\end{equation}
The substitution $q \to q^{-1}$ in \eqref{lambdaeq} gives:
\begin{equation}
 \begin{cases}
   \lambda_1(q^{-1}) + \lambda_2(q^{-1}) + \lambda_3(q^{-1})= a_{[3,1]}(q^{-1})\\
   \lambda_1(q^{-1})\cdot \lambda_2(q^{-1})\cdot \lambda_3(q^{-1}) = (-q)^{-W}\\
   \lambda_1^{-1}(q^{-1}) + \lambda_2^{-1}(q^{-1}) + \lambda_3^{-1}(q^{-1}) = a_{[3,1]}(q)
 \end{cases}
\end{equation}
Let $\lambda_i(q^{-1})\rightarrow \lambda_i^{-1}(q^{-1})=\tilde{\lambda}_i(q)$. Then
\begin{equation}\label{InvLambda}
    \begin{cases}
        \tilde{\lambda}_1 (q) +\tilde{\lambda}_2(q) +\tilde{\lambda}_3(q) =a_{[3,1]}(q)\\
        \tilde{\lambda}_1(q) \cdot\tilde{\lambda}_2(q)\cdot\tilde{\lambda}_3(q)=(-q)^W \\
        \tilde{\lambda}_1^{-1}(q) +\tilde{\lambda}_2^{-1}(q) +\tilde{\lambda}_3^{-1}(q) =a_{[3,1]}(q^{-1})
    \end{cases}
\end{equation}
Comparing~\eqref{InvLambda} with~\eqref{lambdaeq}, we conclude that sets of solutions of these systems coincide:
\begin{equation}\label{TotalLambdaEq}
    \{\lambda_i(q)\}=\{\lambda_i^{-1}(q^{-1})\}\,.
\end{equation}
This can be also proved by consideration of the equivalent characteristic equation on the eigenvalues:
\begin{equation}\label{LambdaCharEq}
    \lambda^3-\lambda^2\cdot a_{[3,1]}(q)+\lambda\cdot q^{-W} a_{[3,1]}(q^{-1})-q^{-W}=0\,.
\end{equation}
Again, after the substitution $q\rightarrow q^{-1}$, $\lambda^{-1}(q^{-1})=\tilde{\lambda}(q)$, one obtains
\begin{equation}
    \tilde{\lambda}^3-\tilde{\lambda}^2\cdot a_{[3,1]}(q)+\tilde{\lambda}\cdot q^{-W} a_{[3,1]}(q^{-1})-q^{-W}=0
\end{equation}
and arrives to~\eqref{TotalLambdaEq}. Equality~\eqref{TotalLambdaEq} means that either $\{\lambda_i(q)\}$ lie on a unit circle, or without loss of generality, $\lambda_1$ lies on a unit circle and $\lambda_2$, $\lambda_3$ lie on a line which intersects the origin. 

If the second case takes place at least for one particular value $k$ from the admissible region~\eqref{4-str-restr}, then the functions $a_{[2,1,1]}$  and $a_{[3,1]}$  are guaranteed not to be constructed from the braid. On the other hand, if all eigenvalues lie on the unit circle, it is a serious argument that the braid does exist.

\begin{figure}[h!]
\scalebox{0.8}{
\begin{picture}(160,80)(-190,-15)
\put(-10,48){\line(1,0){20}}
\put(-10,24){\line(1,0){20}}
\put(-10,0){\line(1,0){80}}
\put(115,0){\line(1,0){75}}

\put(70,-10){\line(1,0){45}}
\put(70,34){\line(1,0){45}}
\put(70,-10){\line(0,1){44}}
\put(115,-10){\line(0,1){44}}
\put(88,9){\mbox{{\large $b_1$}}}

\put(10,14){\line(1,0){45}}
\put(10,58){\line(1,0){45}}
\put(10,14){\line(0,1){44}}
\put(55,14){\line(0,1){44}}
\put(28,33){\mbox{{\large $c_1$}}}

\put(130,14){\line(1,0){45}}
\put(130,58){\line(1,0){45}}
\put(130,14){\line(0,1){44}}
\put(175,14){\line(0,1){44}}
\put(148,33){\mbox{{\large $c_2$}}}

\put(190,-10){\line(1,0){45}}
\put(190,34){\line(1,0){45}}
\put(190,-10){\line(0,1){44}}
\put(235,-10){\line(0,1){44}}
\put(208,9){\mbox{{\large $b_2$}}}

\put(55,24){\line(1,0){15}}
\put(55,48){\line(1,0){75}}

\put(115,24){\line(1,0){15}}

\put(175,24){\line(1,0){15}}

\put(175,48){\line(1,0){80}}

\put(235,24){\line(1,0){20}}
\put(235,0){\line(1,0){20}}

\end{picture}
}
\caption{{\small Considered 3-strand braids which parameterise character expansion coefficients $a_Y$, $-7\leq c_1\,,\,b_1\,,c_2\,,\,b_2\leq 7$ and the withe number $W=2m+5$, $m=0,\dots,5$.}}
\label{fig:four-blocks}
\end{figure}
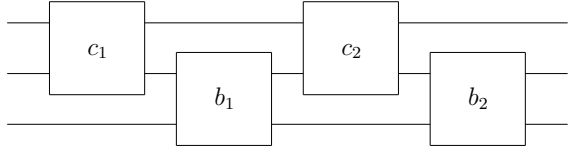

We have considered 3-strand knots of the type $(c_1,b_1;c_2,b_2)$, $-7\leq c_1,b_1,c_2,b_2 \leq 7$, with the writhe number $W=2m+5$, $m = 0, \dots, 5$, see Fig.~\ref{fig:four-blocks}. Recall that these 3-strand braids give values of the coefficient $a_{[2,2]}$. We have found absolute values of eigenvalues $\lambda_i^{[3,1]}(q), \ q=\exp{\frac{2\pi i}{k}}$ for the real values of $k$ from $0.1$ to $150$ in increments of $0.1$ numerically with the help of Wolfram Mathematica. It turns out that the condition $|\lambda^{[3,1]}_i|=1$ at the regions $k \geq 8, \, \text{or}\, k \leq -8, \, \text{or}\, \frac{8}{1+4 z}\leq k\leq \frac{8}{-1+4 z}\, \forall\, z\in \mathbb{Z}/\{0\}$ is fulfilled only for the following braids: $(-5,2;6,2),\, (-2,-1;1,7),\,(-2,-1;3,5),\, (-2,2;3,2),\,(-1,1;2,3),\,(-2,2;1,4)$. Note that the writhe number $W$ of these braids is 5, for higher values of $W$ there are no such braids from the considered set. The closure of each of these braids is a 2-link on 3-strands. Thus, hypothetically, each of them can correspond to a knot on 4 strands, for which the HOMFLY-PT polynomial coincides with~\eqref{difF} at $m=0$.


\subsection{Checks for related knots}\label{sec:rel-knots}
If knots ${\cal K}_m=\overline{\beta}_m$ for the obtained HOMFLY-PT polynomials~\eqref{difF},~\eqref{CycF} do exist, then knots corresponding to $\beta_m^{2k+1}$ and $\beta_m \Delta^{2n}$ with the fundamental half-twist braid $\Delta=\sigma_1 \sigma_2 \sigma_3 \sigma_1 \sigma_2 \sigma_1$ must exist too. In this section, we use the same methods to analyse the HOMFLY-PT invariants for these related knots. 

\subsubsection{Analysis of the HOMFLY-PT polynomials for $\tilde{{\cal K}}_m=\overline{\beta_m \Delta^{2n}}$}
Multiplication of a four-strand braid by the fundamental half-twist braid $\Delta$ raised to an even power changes coefficients of the character expansion of the HOMFLY-PT polynomial in a clear and simple way. Note that $\mathfrak{R}$-matrix representation of $\Delta^{2n}$ in the basis of irreducible representation is
\begin{equation}\label{DeltaProp}
    \left(\Delta^{[2,2]}_{\mathfrak{R}}\right)^{2n}=\mathds{1}\,,\quad \left(\Delta^{[3,1]}_{\mathfrak{R}}\right)^{2n}=q^{4n}\cdot \mathds{1}\,,\quad \left(\Delta^{[2,1,1]}_{\mathfrak{R}}\right)^{2n}=q^{-4n}\cdot \mathds{1}\,.
\end{equation}
Thus, if the $\mathfrak{R}$-matrix representation of a braid $\beta$ is denoted by $\cal R$, then the irreducible blocks of $\tilde{\cal R}={\cal R}\cdot \Delta_{\mathfrak{R}}^{2n}$ are
\begin{equation}\label{R-mat-rep-tilde-beta}
    \tilde{\cal R}^{[2,2]}={\cal R}^{[2,2]}\,,\quad \tilde{\cal R}^{[3,1]}=q^{4n}\cdot {\cal R}^{[3,1]}\,,\quad \tilde{\cal R}^{[2,1,1]}=q^{-4n}\cdot {\cal R}^{[2,1,1]}\,.
\end{equation}
The character expansion of the HOMFLY-PT polynomial for a knot $\tilde{\cal K}=\overline{\beta \Delta^{2n}}$~\eqref{HOMFLYhowR4} gets the form
\begin{equation}\label{HDelta}
    H^{\tilde{{\cal K}}}= A^{-12n}\cdot A^{-W}\left(q^{12n}\cdot q^W S^{*}_{[4]} + q^{4n}\cdot a_{[3,1]}S^{*}_{[3,1]}+ a_{[2,2]}S^{*}_{[2,2]}+ q^{-4n}\cdot a_{[2,1,1]}S^{*}_{[2,1,1]} - q^{-12n}\cdot q^{-W} S^{*}_{[1,1,1,1]}\right)\,,
\end{equation}
where we assume that the writhe number $W$ of a braid $\beta$ is odd, and
\begin{equation}
    a_{[2,2]}=\tr \left({\cal R}^{[2,2]}\right)\,,\quad a_{[3,1]}=\tr\left( {\cal R}^{[3,1]}\right)\,,\quad a_{[2,1,1]}=\tr\left( {\cal R}^{[2,1,1]}\right)\,.
\end{equation}

\subsubsection*{Eigenvalues}
The first idea is to check whether the eigenvalues of $\mathfrak{R}$-matrix representation of our hypothetical knots $\tilde{{\cal K}}_m=\overline{\beta_m \Delta^{2n}}$ lie on a unit circle at the point $q=\exp \frac{2\pi i}{k}$, $k\geq 8$ and $k=1,2$. However, this check is trivial. After multiplication by $\Delta^{2n}$, eigenvalues of $\mathfrak{R}$-matrix of $\beta_m$ change only by a phase factor according to~\eqref{R-mat-rep-tilde-beta}. Thus, if eigenvalues of $\mathfrak{R}$-matrix of $\beta_m$ lie on a unit circle at the special roots of unity, then the eigenvalues of $\mathfrak{R}$-matrix of $\overline{\beta_m \Delta^{2n}}$ also lie on a unit circle.

\subsubsection*{Differential expansion}
The second idea is to check whether the HOMFLY-PT polynomial~\eqref{HDelta} with coefficients $a_Y$ defined by~\eqref{a-coef} possesses the structure of differential expansion~\eqref{difexp2}. Indeed, by straightforward calculation, one can get sure that
\begin{equation}
    H^{\tilde{{\cal K}}_m}(q,A=q)-1=0 \quad \text{and} \quad H^{\tilde{{\cal K}}_m}(q^{-1},A=q)-1=0\,.
\end{equation}

\subsubsection*{Perturbative expansion}
Having proved the validness of the differential expansion for $H^{\tilde{{\cal K}}_m}$, to show that these HOMFLY-PT polynomials can be perturbative expanded, one should check that $H^{\tilde{{\cal K}}_m}(q,A)=H^{\tilde{{\cal K}}_m}(q^{-1},A)$ (see Section~\ref{sec:PertExpDef}). After inversion of $q$ the HOMFLY-PT polynomial~\eqref{HDelta} transforms as follows:
{\small \begin{equation}
   H^{\tilde{{\cal K}}}(q^{-1},A) = -A^{-W-12n}\left(q^{-W-12n}S^*_{[1,1,1,1]}+a_{[3,1]}(q^{-1})q^{-4n}S^*_{[2,1,1]}+a_{[2,2]}(q^{-1})S^*_{[2,2]}-q^{12n+W}S^*_{[4]}+a_{[2,1,1]}(q^{-1})q^{4n}S^*_{[3,1]}\right)\,.
\end{equation}}
Taking into account condition~\eqref{a-rel}, we get that the equality $H^{\tilde{{\cal K}}_m}(q^{-1},A)=H^{\tilde{{\cal K}}_m}(q,A)$ is automatically fulfilled, and thus, does not give additional constraints on $a_Y$-coefficients~\eqref{a-coef}.

\subsubsection{Analysis of the HOMFLY-PT polynomials for ${\cal K}^{2k+1}=\overline{\beta^{2k+1}}$}
Given the HOMFLY-PT polynomial for a knot $\cal K=\overline{\beta}$ in the form of the character expansion~\eqref{Hforcoef}, one can construct the HOMFLY-PT polynomial for a knot ${\cal K}^{2k+1}=\overline{\beta^{2k+1}}$:
{\small \begin{equation}\label{HOMFLYK2k+1}
\begin{aligned}
    H^{{\cal K}^{2k+1}}&=A^{-(2k+1)W}\Big\{q^{(2k+1)W}S^*_{[4]}+(-q^{-1})^{(2k+1)W}S^*_{[1,1,1,1]}+\left(\left(\lambda_1^{[3,1]}\right)^{2k+1}+\left(\lambda_2^{[3,1]}\right)^{2k+1}+\left(\lambda_3^{[3,1]}\right)^{2k+1}\right)S^*_{[3,1]}+\\
    &+\left(\left(\lambda_1^{[2,2]}\right)^{2k+1}+\left(\lambda_2^{[2,2]}\right)^{2k+1}\right)S^*_{[2,2]}+\left(\left(\lambda_1^{[2,1,1]}\right)^{2k+1}+\left(\lambda_2^{[2,1,1]}\right)^{2k+1}+\left(\lambda_3^{[2,1,1]}\right)^{2k+1}\right)S^*_{[2,1,1]}\Big\}\,,
\end{aligned}
\end{equation}}
where we assume that the writhe number $W$ of a knot $\cal K$ is odd, and $\lambda_i^Y$ are eigenvalues of $\mathfrak{R}$-matrix representation of a braid $\beta$. These eigenvalues can be derived from systems~\eqref{lambda4} and~\eqref{lambda22-sys} and within the condition $\lambda^{[3,1]}_i(q)=-\lambda^{[2,1,1]}_j(q^{-1})$.

\subsubsection*{Eigenvalues}
If eigenvalues of $\mathfrak{R}$-matrix representation ${\cal R}$ of a braid $\beta$ lie on a unit circle, then it is obvious that eigenvalues of ${\cal R}^{2k+1}$ also lie on a unit circle. Thus, this check is trivial.

\subsubsection*{Differential expansion} 
One can show by a straightforward calculation that for our coefficients $a_Y$~\eqref{a-coef}
\begin{equation}\label{HK2k+1DE}
    H^{{\cal K}^{2k+1}_m}=1+(A^{2} - q^{2})(A^{2} - q^{-2})\cdot C^{{\cal K}^{2k+1}_m}\,.
\end{equation}

\subsubsection*{Perturbative expansion}
Having~\eqref{HK2k+1DE}, it is enough to show that $H^{{\cal K}^{2k+1}_m}(q, A)=H^{{\cal K}^{2k+1}_m}(q^{-1},A)$ to prove that the HOMFLY-PT polynomial $H^{{\cal K}^{2k+1}}(q, A)$ has right perturbative expansion. Change $q\rightarrow q^{-1}$ transorms the HOMFLY-PT polynomial~\eqref{HOMFLYK2k+1} to
\begin{equation}
\begin{aligned}
    H^{{\cal K}^{2k+1}}(q^{-1},A)&=-A^{-(2k+1)(2m+5)}\Big\{q^{-(2k+1)(2m+5)}S^*_{[1,1,1,1]}+(-q)^{(2k+1)(2m+5)}S^*_{[4]}+\\ &+\left(\left(\lambda_1^{[3,1]}(q^{-1})\right)^{2k+1}+\left(\lambda_2^{[3,1]}(q^{-1})\right)^{2k+1}+\left(\lambda_3^{[3,1]}(q^{-1})\right)^{2k+1}\right)S^*_{[2,1,1]}+\\
    &+\left(\left(\lambda_1^{[2,2]}(q^{-1})\right)^{2k+1}+\left(\lambda_2^{[2,2]}(q^{-1})\right)^{2k+1}\right)S^*_{[2,2]}+\\
    &+\left(\left(\lambda_1^{[2,1,1]}(q^{-1})\right)^{2k+1}+\left(\lambda_2^{[2,1,1]}(q^{-1})\right)^{2k+1}+\left(\lambda_3^{[2,1,1]}(q^{-1})\right)^{2k+1}\right)S^*_{[3,1]}\Big\}\,.
\end{aligned}
\end{equation}
Baring in mind relation~\eqref{LambdaRel}, we see that the HOMFLY-PT polynomial $H^{{\cal K}^{2k+1}_m}$ is invariant under inversion of $q$, and we get no additional constraints.

\subsection{Positive braids}\label{sec:pos-braids}
So, as we have seen, the obtained formulas for the hypothetical HOMFLY-PT polynomials~\eqref{difF},~\eqref{CycF} confidently withstand all checks. This inspires great optimism and raises the question of how to find such a braid and knot, respectively. We conclude this Section with a discussion of this issue.

We want to find any 4-strand braid that, in the representation [3,1], satisfies conditions \eqref{a-coef}. For simplicity, we consider here only $W=5$. In order to efficiently enumerate only braids with fixed writhe number, we do the following. 

Recall that any braid $\beta$ can be represented in the form (formula),
\begin{equation}\label{betabraid}
\beta = \beta_{+}\Delta^{-p}\,, 
\end{equation}
where $\beta_+$ is some positive braid, $\Delta$ is the fundamental half-twist braid (which generates the center of $\mathcal{B}_n$) and $p \in \mathbb{Z}$. 

Let us consider all positive braids $\beta_+ \in B_4$ with the writhe number $W = 5 + 6p$. Since  $\Delta=\sigma_1 \sigma_2 \sigma_3 \, \sigma_1 \sigma_2 \, \sigma_1$, then its writhe number is 6. Therefore, enumerating all positive braids with $W = 5 + 6p$, we enumerate all braids $\beta$ with $W=5$. Further we note that it is enough to consider only even $p$, what gives us two properties. First, the addition of $\Delta^{2}$ does not change the number of components of the initial knot or link, because $\Delta^{2}$ is a pure braid. Second, character expansions \eqref{Hforcoef} of HOMFLY-PT polynomial for the knot $\overline{\beta}$ and for $\overline{\beta_+}$ are related in a simple way due to~\eqref{DeltaProp}:
\begin{equation}\label{Newacoefficient}
a^+_{[2,2]} = a_{[2,2]}\,, \quad a^+_{[3,1]} = q^{4p} \, a_{[3,1]}\,, \quad a^+_{[2,1,1]} = q^{-4p} \, a_{[2,1,1]}\,,
\end{equation}
where $a^+_Y$ are coefficients of the character expansion for the HOMFLY-PT polynomial for the braid $\beta_+=\beta \Delta^p$, and coefficients $a_Y$ are given by~\eqref{a-coef}. 

Thus, the proposed search method is as follows. First, we fix the number of full twists given by $\Delta$ to be equal to $n$, i.e. $p=2k$. This corresponds to positive braids with the number of intersections equal to $5 + 12k$. There is a finite number of such braids, equal to $3^{5+12k}$. Second, we select from them only those that correspond to knots. This can be done by looking at the permutation that the braid naturally induces. Third, for the resulting braids, we calculate the coefficients of the character expansion $a^{+}_Y$ and check relation \eqref{a-coef}, taking into account formulas \eqref{Newacoefficient}. Any braid constructed in this way is suitable for us. The reverse is also true. And if such a braid exists, then we will find it in a finite number of steps. Note that the simplest case $k=1$ is already interesting, because the resulting braid $\beta = \beta_+ \Delta^{-2}$ can contain $29 = 17 + 12$ crossings that exceeds the value in  \cite{tuzun2021verification}.

\section*{Ackhowledgements}

We are grateful for enlightening discussions to V. Mishnyakov, A. Morozov, And. Morozov and N. Tselousov.

This work was funded by the Russian Science Foundation (Grant No.20-71-10073).

\printbibliography

\end{document}